\documentclass[11pt,usenames,dvipsnames]{amsart}

\usepackage{ amssymb }
\usepackage{amsfonts}
\usepackage{tikz}
\usepackage{xcolor}
\usepackage{comment}
\usepackage{mathrsfs}
\usepackage{ifdraft}
\usepackage[colorinlistoftodos,prependcaption,textsize=tiny]{todonotes}

\usepackage{enumerate}
\usepackage{enumitem}
\usepackage{amsmath,amscd}
\usepackage{amsmath,amsthm,amssymb,anysize}
\usepackage[numbers,sort&compress]{natbib}
\usepackage{color}
\usepackage{graphicx}
\usepackage[all,cmtip]{xy}
\usepackage{tikz-cd} 
\usepackage{blindtext}
\usepackage[normalem]{ulem}
\usepackage{hyperref}

\newtheorem{thm}{Theorem}
\newtheorem*{thm*}{Theorem}
\newtheorem{lem}[thm]{Lemma}

\newtheorem{prop}[thm]{Proposition}
\newtheorem{cor}[thm]{Corollary}
\newtheorem{nota}[thm]{Notations}
\newtheorem{conj}[thm]{Conjecture}
\newtheorem{que}[thm]{Question}

\newtheorem{defn}[thm]{Definition}
\newtheorem{remark}[thm]{Remark}
\newcommand{\opn}{\operatorname}
 \newcommand{\Zp}{\mathbb{Z}_{(p)}}
\newcommand{\K}[1]{K(\mathbb{Z}/#1, 2)}

\newcommand{\xx}[2]{x_{#1}^{(#2)}}
\newcommand{\yy}[2]{y_{#1}^{(#2)}}

\newcommand{\xxx}[3]{x_{#1,#2}^{(#3)}}
\newcommand{\yyy}[3]{y_{#1,#2}^{(#3)}}

\newcommand{\xxi}[2]{\xi_{#1, #2}}

\newcommand{\bxxi}[2]{\bar{\xi}_{#1, #2}}
\newcommand{\eeta}[2]{\eta_{#1, #2}}

\newcommand{\zzeta}[2]{\zeta_{#1}^{(#2)}}
\newcommand{\zzzeta}[3]{\zeta_{#1, #2}^{(#3)}}

\newcommand{\barxi}{\bar{\xi}}
\newcommand{\bareta}{\bar{\eta}}
\newcommand{\ssig}[2]{\sigma_{#1}^{(#2)}}
\newcommand{\ttau}[2]{\tau_{#1}^{(#2)}}
\newcommand{\om}[3]{\omega_{#1}^{#2,#3}}
\newcommand{\bom}[3]{\bar{\omega}_{#1}^{#2,#3}}
\newcommand{\uinf}[1]{u_{#1,\infty}}
\newcommand{\vinf}[1]{v_{#1,\infty}}
\newcommand{\mmu}[2]{\mu_{#1}^{(#2)}}

\newcommand{\cupo}{\smile_1}

\newcommand{\per}{\operatorname{per}}
\newcommand{\ind}{\operatorname{ind}}
\newcommand{\Br}{\operatorname{Br}_{top}}
\newcommand{\sk}[2]{\operatorname{sk}_{#1}K(\mathbb{Z}/#2,2)}
\newcommand{\wE}{\widetilde{E}}
\newcommand{\wdd}{\widetilde{d}}
\newcommand{\td}{\widetilde}
\newcommand{\pp}{\operatorname{P}}
\newcommand{\htf}[1]{\widehat{f}_{#1}}

\marginsize{3cm}{3cm}{2.5cm}{2.2cm}
\numberwithin{equation}{section}
\numberwithin{thm}{section}
\def\Z{\mathbb{Z}}

\def\Q{\mathbb{Q}}

\def\C{\mathbb{C}}

\def\O{\mathfrak{O}}

\def\F{\mathbb{F}}

\usepackage{color}

\title{The $p$-local topological period-index problem in dimensions $2p+2$ and below}

\author{Diarmuid Crowley}
\address{School of Mathematics and Statistics
The University of Melbourne
Parkville, VIC, 3010, Australia}
\email{dcrowley@unimelb.edu.au}

\author{Xing Gu}
\address{Institute for Theoretical Sciences, School of Sicence, Westlake University, 600 Dunyu Road, Xihu District, Hangzhou, Zhejiang 310030, China}
\email{guxing@westlake.edu.cn}

\address{}
\email{}

\date{\today}

\begin{document}
\begin{abstract}
    For any odd prime number 
    $p$, we solve the topological period-index problem for a $p$-primary element of the Brauer group of a CW-complex of dimension not greater than $2p{+}2$.
\end{abstract}
\maketitle
\vskip 10pt

\tableofcontents

\section{Introduction}\label{sec:intro}
Let $X$ be a topological space. By the \textit{topological Brauer} group $X$, denoted by $\Br(X)$, we mean the subgroup $TH^3(X;\Z) \subseteq H^3(X; \Z)$ of all torsion elements.

Equivalently, $\Br(X)$ may be identified as the group whose underlying set is the Azumaya algebras (i.e., bundles of matrix algebras over $\C$) modulo the Brauer equivalence: $\mathscr{A}_0$ and $\mathscr{A}_1$ are called Brauer equivalent if there are vector bundles $\mathscr{E}_0$ and $\mathscr{E}_1$ such that
$$\mathscr{A}_0\otimes\operatorname{End}(\mathscr{E}_0)\cong \mathscr{A}_1\otimes\operatorname{End}(\mathscr{E}_1),$$
such that tensor product gives the abelian group structure on $\Br(X)$. 
We briefly explain why the above definitions coincide. Indeed, 
Azumaya algebras over $X$ of degree $r$ are classified by the collection of $PU(r)$-torsors over $X$, i.e., the cohomology set $H^{1}(X;PU(r))$, where $PU(r)$ is the projective unitary group of degree $r$. Consider the short exact sequences of Lie groups
\begin{equation}\label{ses1}
1\rightarrow S^1\rightarrow U(r)\rightarrow PU(r)\rightarrow 1
\end{equation}
and
\begin{equation}\label{ses2}
0\rightarrow\mathbb{Z}\rightarrow\mathbb{C}\xrightarrow{\textrm{exp}} S^1\rightarrow 1
\end{equation}
where the arrow $S^1\rightarrow U_r$ is the inclusion of scalars. Then the composition of Bockstein homomorphisms
\begin{equation}\label{Bockstein}
H^{1}(X;PU(r))\rightarrow H^{2}(X;S^1)\rightarrow H^{3}(X;\mathbb{Z})
\end{equation}
associates an Azumaya algebra $\mathscr{A}$ to a class $\alpha\in H^{3}(X;\mathbb{Z})$. 

Serre showed (\cite{Gr}) that when $X$ is a finite CW-complex, for any $\alpha\in H^{3}(X;\mathbb{Z})_{\textrm{tor}}$, there is some $r$ such that a $PU(r)$-torsor over $X$ is associated to $\alpha$ via the homomorphism (\ref{Bockstein}). Let $\operatorname{per}(\alpha)$ denote the order of $\alpha$ as an element of the group $H^{3}(X;\mathbb{Z})$, then Serre also showed (\cite{Gr}) that $\operatorname{per}(\alpha)|r$, for all $r$ such that there is a $PU(r)$-torsor over $X$ associated to $\alpha$ in the way described above. Let $\operatorname{ind}(\alpha)$ denote the greatest common divisor of all such $r$, then in particular we have
\begin{equation}\label{per div ind}
\operatorname{per}(\alpha)|\operatorname{ind}(\alpha).
\end{equation}
Furthermore, Antieau and Williams showed  the following
\begin{thm}[Theorem 6, \cite{antieau2015prime}]\label{same prime divisors}
Let $(X,R)$ be a locally-ringed connected topos and let $\alpha\in\operatorname{Br}(X,R)$. There exists a representative $A$ of $\alpha$ such that the prime numbers dividing $\operatorname{per}(\alpha)$ and $\operatorname{deg}(A)$ coincide.
\end{thm}
In particular, the period and index of any $\alpha\in\operatorname{Br}(X)$ have the same prime divisors, for a connected finite CW-complex $X$. Hence, for a sufficiently large integer $e$ we have
\begin{equation}\label{ind div per}
\operatorname{ind}(\alpha)|\operatorname{per}(\alpha)^e.
\end{equation}
One version of the topological period-index problem can be stated as follows:

\begin{que}\label{que:tpip}
    For (connected, finite) CW-complexes of a given dimension $d$, find the sharp lower bound of $e$ such that (\ref{ind div per}) holds for all finite CW-complex $X$ in $\textrm{C}$ and all elements $\alpha\in\Br(X)$.
\end{que}

An upper bound of the index in terms of the period and the dimension of the ambient CW-complex has been established by Antieau and Williams as follows:
\begin{thm}[Theorem A, \cite{antieau2021topological}]\label{thm:big_upper_bound}
    Let $X$ be a $2d$-dimensional CW-complex, and let $\alpha\in\Br(X)$. Setting $n = \per(\alpha)$, we have
    \[\ind(\alpha)\mid n^{d-1}\prod_{p\mid n }p^{\nu_p((d-1)!)}.\]
\end{thm}

The topological period-index problem is motivated by its analog in algebraic geometry, where the topological Brauer group of a CW-complex is replaced by the usual Brauer group of a scheme, and where the period and index are defined similarly. We have the following folklore conjecture:
\begin{conj}[Colliot-Th{\'e}l{\`e}ne]\label{algconj}
Let $k$ be either a $C_d$-field or the function field of a $d$-dimensional variety over an algebraically closed field. Let $\alpha\in \textrm{Br}(k)$, and suppose that $\operatorname{per}(\alpha)$ is prime to the characteristic of $k$. Then
$$\operatorname{ind}(\alpha)|\operatorname{per}(\alpha)^{d-1}.$$
\end{conj}

Our first main result, Theorem~\ref{thm:cw} below, says that the upper bound established in Theorem~\ref{thm:big_upper_bound} is sharp, when the least prime divisor of the period is large relative to $d$.

\begin{thm}\label{thm:cw}
    Let $p$ be an odd prime number. Let $d$ be an integer such that $1\leq d\leq 2p+2$.
    \begin{enumerate}
        \item Let $X$ be a CW-complex of dimension $d$.  Let $\alpha\in\Br(X)$ be such that the least prime divisor of $\per(\alpha)$ is no less than $p$. Then we have
        \begin{equation*}
            \begin{cases}
                \ind(\alpha)\mid \per(\alpha)^{[(d+1)/2]-1},\ d < 2p + 1,\\
                \ind(\alpha)\mid p\per(\alpha)^{[(d+1)/2]-1},\ d = 2p +1 \textrm{ or }2p+2.
            \end{cases}
        \end{equation*}
        
        \item For $n$ a natural number such that its least prime divisor is no less than $p$, there is a CW-complex $X$ of dimension $d$, and $\alpha\in\Br(X)$ such that $\per(\alpha) = n$ and 
        \begin{equation*}
            \begin{cases}
                \ind(\alpha) = \per(\alpha)^{[(d+1)/2]-1},\ d < 2p + 1,\\
                \ind(\alpha) = p\per(\alpha)^{[(d+1)/2]-1},\ d = 2p +1 \textrm{ or }2p+2.
            \end{cases}
        \end{equation*}
    \end{enumerate}  
\end{thm}
\begin{remark}\label{rem:d}
   One should think of the exponent $[(d+1)/2]-1$ as playing the role of the ``complex dimension'' of $X$, in the sense that when $X$ is a complex manifold or a complex variety, $d$ is even and $[(d+1)/2] = d$. This is consistent with our choice of working with $PU(r)$-principal bundles. 
\end{remark}

Our proof of Theorem~\ref{thm:cw} does not apply to the case $p = 2$, but it has been  solved by Antieau-Williams \cite{antieau2014topological}. 
The case $p=3$ is solved by Gu \cite{gu2019topological}.

In Question~\ref{que:tpip}, if we replace ``(connected, finite) CW complexes'' with ``closed, smooth, oriented manifold'', we may improve the index bound for $d = 2p+2$, so that in the case $d\leq 2p+2$, and $d$ is even,  the index bound is consistent with Conjecture\ref{algconj}, in the sense of Remark~\ref{rem:d}.

\begin{thm}\label{thm:manifold}
    Let $p$ be an odd number. Let $M$ be a closed oriented 
    smooth
    manifold of dimension $d$ for $0\leq d\leq 2p+2$. Let $\alpha$ be a topological Brauer class of $M$ such that the least prime divisor of $\per{\alpha}$ is no less than $p$. Then we have 
    \[\ind(\alpha)\mid\per(\alpha)^{[(d+1)/2]-1}.\]
\end{thm}

As formualated by Antieau and Williams \cite{antieau2014period},
the topological period-index problem can be viewed as a problem in twisted $K$-theory.  For $X$ a finite $CW$-complex and a``twist" $\alpha \in H^3(X; \Z)$,
there is a twisted Atiyah-Hirzebruch spectral sequence $(E,d)$ converging to 
$KU^*_{\alpha}(X)$, the $\alpha$-twisted complex $K$-theory of $X$.
As we recall in Section~\ref{sec:twistedAH}, Antieau and Williams \cite{antieau2014period} showed that the period and index can be interpreted in terms of this spectral sequence.  
The next theorem is essentially another way of stating Theorem~\ref{thm:cw}.

\begin{thm}\label{cor:AHSS}
    Let $\alpha\in H^{3}(K(\Z/n,2);\Z)$ be the Bockstein image of the fundamental class. Then in the twisted Atiyah-Hirzebruch spectral sequence for $KU^0_{\alpha}(K(\Z/n,2))$, the differentials
    \[d_{2r+1}^{0,0}\colon E_{2r+1}^{0,0}\to E_{2r+1}^{2r+1, -2r} = E_2^{2r+1,-2r}\cong H^{2r-1}(K(\Z/n,2);\Z)\]
    are surjective for $r\leq p$.
\end{thm}

Yet another way of (partially) interpreting Theorem~\ref{thm:cw} is in terms of the cohomology of the classifying spaces of some Lie groups. Notice that the center of $SU(m)$ is the cyclic group $\mu_m$. For $n\mid m$, let $\mu_n$ be the subgroup of $\mu_m$ generated by $mn^{-1}$. Let $G(n,m) = SU(m)/\mu_n$. 
The short exact sequence 
\[1\to\mu_n\to SU(m)\to G(n,m)\to 1\]
induces a map $BG(n,m)\xrightarrow{\alpha_{n,m}}K(\Z/n,2)$.
\begin{thm}\label{thm:Gnm} 
    Suppose one of the following is true:
    \begin{enumerate}
        \item $r < p$, $n$ a power of $p$, $n^r\mid m$, or
        \item $r= p$, $n$ a power of $p$, $pn^r\mid m$. 
    \end{enumerate}
    Then the map
    \[\alpha_{n,m}^*\colon H^{\leq 2r+2}(K(\Z/n,2);\Z)\to  H^{\leq 2r+2}(BG(n,m);\Z)\]
    is injective.
\end{thm}

The key topological construction of this work is as follows. Consider the lifting diagrams of the form
\begin{equation}\label{eq:lift_diag_general}
    \begin{tikzcd}
       & BPU(m)\arrow[d]\\
    X\arrow[r,"\alpha"]\arrow[ur,dashed]& K(\Z,3)      
    \end{tikzcd}
\end{equation}
for various $m$, where $\alpha$ represents an element in $\Br(X)$, and the vertical arrow represents a generator of $\Br(BPU(m))$. The key idea involved in the proof of Theorem~\ref{thm:cw} is to study the cohomological obstructions to the liftings indicated by the dashed arrow in \eqref{eq:lift_diag_general} via Massey higher products. 

The rest of this paper is organised as follows.

In Section~\ref{sec:twistedAH}, we recall the Atiyah-Hirzebruch spectral sequence for twisted complex topological K-theory and its connection to the topological period-index problem.

In Section~\ref{sec:Massey}, we briefly review the theory of Massey higher products following \cite{kraines1966massey}.  

The ``universal examples'' of the CW-complex $X$ in \eqref{eq:lift_diag_general} is the $2d$-skeletons of the Eilenberg-Mac Lane spaces $K(\Z/n,2)$. We review their cohomology in Section~\ref{sec:Kn2}.

Another key ingredient is the cohomology of $BPU(m)$, which is the subject of many works including \cite{gu2019cohomology}, \cite{gu2021p}, \cite{crowley2026}, \cite{fan2025formula}. We revisit this subject in Section~\ref{sec:cohomology}.

With the aforementioned preliminaries established, we carry out the proof of Theorem~\ref{thm:cw} in Section~\ref{sec:existence}.

In Section~\ref{sec:Steenrod}, \ref{sec:linking}, \ref{sec:bilinear} and \ref{sec:2p+2}, we prove Theorem~\ref{thm:manifold}, which says that the topological period-index problem tends to have a better answer for closed oriented differentiable manifolds then for CW-complexes. 

In Section~\ref{sec:G}, we prove Theorem~\ref{thm:Gnm}, which is a consequence of Theorem~\ref{thm:cw}.

{\bf Acknowledgements: }  The authors would like to thank Ben Antieau, Ben Williams, Christian Haesemeyer, Yang Su, Songyan Xie, and Foling Zou for helpful conversations.

{\bf Remark: }
The proofs of Theorems~\ref{thm:cw} and
\ref{thm:manifold} were completed in 2025.
Subsequently, 
we spent some time on an unsuccessful attempt to construct closed oriented smooth manifolds realising the index bound of 
Theorem~\ref{thm:manifold}.

{\bf AI statement: } No AI assistance was used in writing this paper.

\section{The twisted Atiyah-Hirzebruch spectral sequence}\label{sec:twistedAH}
In this section we recall some facts on the Atiyah-Hirzebruch spectral sequence for twisted K-theory, and we interpret period and index in terms of invariants of this spectral sequence.

For an integer $t$, let 
\begin{equation*}
    \Z(t/2) = 
    \begin{cases}
        \Z,\ t\textrm{ even},\\
        0,\ t\textrm{ odd}.
    \end{cases}
\end{equation*}
\begin{prop}[\cite{antieau2014period}, Proposition~2.3, Proposition~2.4]\label{pro:per_ind_ahss}
    Let $X$ be a topological space and let $\alpha\in\Br(X)$. Then there is a spectral sequence $(E,d)$ with
    \[E_2^{s,t}\cong H^s(X;\Z(t/2))\Rightarrow KU^{s+t}_{\alpha}(X)\]
    and
    \[d_r^{s,t}\colon E_r^{s,t}\to E_r^{s+r,t+r-1},\ d_3(1) = \alpha,\]
    which converges strongly if $X$ is a finite CW-complex, and we have
    \[E_{\infty}^{0,0} = \ind(\alpha)E_2^{0,0}.\] 
\end{prop}

The twisted Atiyah-Hirzebruch spectral sequence has a bi-graded module structure over the (un-twisted) Atiyah-Hirzebruch spectral sequence for $K_{\alpha}(X)$, denoted by $(\widetilde{E}, \widetilde{d})$ and we have

\begin{prop}[Product formula]\label{pro:product}
    For $a\in\widetilde{E}_r^{s,t}, b\in E_r^{s',t'}$, we have
    \[d_r^{s+s',t+t'}(ab) = \widetilde{d}_r^{s,t}(a)b + (-1)^{\deg(a)}a d_r^{s',t'}(b).\]
\end{prop}

\section{Massey higher products}\label{sec:Massey}
We recollect the notion of Massey higher products following \cite{kraines1966massey}. Throughout this section, we fix a topological space $X$ and a commutative unital ring $R$, and let $d$ denote the differntial of the cochain complex $C^*(X;R)$,  regarded as a differential graded algebra. For $\alpha\in C^k(X;R)$, we write $\bar{\alpha}$ for $(-1)^{\deg{\alpha}}$.
\begin{defn}\label{def:defining_system}
   Let $\alpha$, $1\leq i\leq r+1$ be a collection of homogeneous cocycles in $C^*(X;R)$. Recall that a \textit{defining system for Massey products} of $\alpha$, $1\leq i\leq r+1$ is a collection of homogeneous elements in $C^*(X;R)$, written as
    \begin{equation}\label{eq:def_sys}
        \begin{matrix}
            \alpha_1 & \beta_{1,1} & \beta_{1,2} & \cdots &\cdots &*\\
               & \alpha_2   & \beta_{2,2} & \cdots &\cdots &\beta_{2,r}\\
               &          & \alpha_3   & \cdots &\cdots &\beta_{3,r}\\
               &          &          & \ddots &\cdots &\vdots\\
               &          &          &        &\alpha_r &\beta_{r,r}\\
               &          &          &        &       &\alpha_{r+1}  
        \end{matrix}      
    \end{equation}
    such that for $i<j$, we have
    \[d(\beta_{i,j}) = \bar{\alpha}_i\beta_{i+1,j} + \sum_{k=i}^{j-2}\bar{\beta}_{i,k}\beta_{k+1,j} + \bar{\beta}_{i,j-1}\alpha_j,\]
    and 
    \[d(\beta_{i,i}) = \bar{\alpha_i}\alpha_{i+1}.\]
    The \textit{Massey product} of $[\alpha_1],\cdots,[\alpha_{i+1}]$ is the subset of $H^*(X,R)$ defined by 
    \[\langle[\alpha_1],\cdots,[\alpha_{r+1}]\rangle = \{[\bar{\alpha}_1\beta_{2,r} + \sum_{j=2}^{r-2}\bar{\beta}_{1,j}\beta_{j+1,r} + \bar{\beta}_{1,r}\alpha_{r+1}]\mid \alpha_i,\beta_{i,j}\textrm{ form a defining system.}\}\]
    for $r\geq 2$ and 
    \[\langle[\alpha_1],[\alpha_2]\rangle = [\alpha_1][\alpha_2].\]
\end{defn}
\begin{defn}\label{def:acyclic}
   Let $\alpha$, $1\leq i\leq r+1$ be a collection of cocycles in $A$, such that the Massey higher product
    \[M = \langle[\alpha_1],\cdots,[\alpha_{r+1}]\rangle\]
   is defined, and suppose there is a defining system \eqref{eq:def_sys}. Then an \textit{acyclic defining system} for $M$ is of the form 
     \begin{equation*}
        \begin{matrix}
            \alpha_1 & \beta_{1,1} & \beta_{1,2} & \cdots &\cdots &\beta_{1,r}\\
               & \alpha_2   & \beta_{2,2} & \cdots &\cdots &\beta_{2,r}\\
               &          & \alpha_3   & \cdots &\cdots &\beta_{3,r}\\
               &          &          & \ddots &\cdots &\vdots\\
               &          &          &        &\alpha_r &\beta_{r,r}\\
               &          &          &        &       &\alpha_{r+1}  
        \end{matrix}      
    \end{equation*}
    such that  
    \[d(\beta_{1,r}) =\bar{\alpha}_1\beta_{2,r} + \sum_{j=2}^{r-2}\bar{\beta}_{1,j}\beta_{j+1,r} + \bar{\beta}_{1,r}\alpha_{r+1}. \]
\end{defn}
The existence of an acyclic defining system for $M$ means, in particular, that $0\in M$. 
\begin{defn}\label{def:normal}
    Let $a\in H^*(A,d)$ be a homogeneous class represented by a cocycle $\alpha\in A$. For $r\geq 1$, an acyclic defining system for the $r$-fold Massey product $\langle a\cdots,a\rangle$ of the form 
    \begin{equation}\label{eq:nor_def_sys}
        \begin{matrix}
            \alpha_1 & \alpha_2 & \alpha_3 & \cdots &\cdots &\alpha_{r}\\
               & \alpha_1   & \alpha_2 & \cdots &\cdots &\alpha_{r-1}\\
               &          & \alpha_1   & \cdots &\cdots &\alpha_{r-2}\\
               &          &          & \ddots &\cdots &\vdots\\
               &          &          &        &\alpha_1 &\alpha_2\\
               &          &          &        &       &\alpha_1      
        \end{matrix}      
    \end{equation}
    where $\alpha_1 = \alpha$ is called an $r$-fold normal acyclic defining system for $[\alpha]$. When the context is clear, we suppress $r$ and say normal acyclic defining systems.    
\end{defn}

The following is a mild generalization of Theorem~15 of \cite{kraines1966massey}.
\begin{lem}\label{lem:nor_def_sys_ext}
    Let $r > 0$ be an integer. Suppose that $R$ contains $\frac{1}{k}$ for $0 < k \leq r$. 
    Let $a\in H^*(X;R)$ be an odd-degree homogeneous class represented by a cocycle $\alpha\in C^*(X;R)$. Then there are  $\alpha_{k}\in C^*(X;R)$, $0 < k\leq r$ such that we have an    
   $r$-fold normal acyclic defining system for $x$ of the form
         \begin{equation*}
        \begin{matrix}
             \alpha_1 & \alpha_2 & \cdots &\cdots &\alpha_{r-1} &\alpha_r\\
                & \alpha_1 & \cdots &\cdots &\alpha_{r-2}     &\alpha_{r-1}\\
               &         & \alpha_1   & \alpha_1 &\cdots &\alpha_{r-3}    \\
               &          &          & \ddots &\cdots &\vdots           \\
               &          &          &        &\alpha_1 &\alpha_2          \\
               &          &          &        &       &\alpha_1           \\
               &          &          &        &       &                 
        \end{matrix}      
    \end{equation*}
    where $\alpha_1 = \alpha$.
\end{lem}
\begin{proof}
   The proof of Theorem~15 of \cite{kraines1966massey} applies verbatim, and we briefly recall it for the sake of completeness. Inductively, let $\alpha_i = \frac{1}{i}\alpha_{i-1}\cupo \alpha_1$.
    Recall (2.5) of \cite{kraines1966massey}, the differential for cup-$1$ products 
    \[d(x\cupo y) = \bar{x}\bar{y} + (-1)^{\deg x \deg y}\bar{y}\bar{x} + d(x)y + \bar{x}d(y),\]
    and (2.6) of \cite{kraines1966massey}, Hirsch's formula
    \[(xy)\cupo z = x(y\cupo z) + (-1)^{\deg y(\deg z -1)}(x\cupo y)z.\]  
    An inductive argument then shows
    \[d(\alpha_r) = -\sum_{i=1}^{r-1}\alpha_i\alpha_{r-i}\]
    and the conclusion follows.
\end{proof}


\section{The cohomology of $\K{n}$}\label{sec:Kn2}
For an odd prime $p$, let $n=p^e$ for some positive integers $e$.
 

Throughout the rest of this paper, we use the following notations for various operations on cohomology. 
\begin{nota}\label{not:operations}
    \begin{itemize}
        \item Let $\rho_k\colon H^*(\square;\Zp)\to H^*(\square;\Z/k)$ denote the mod $k$-reduction.
        \item For $k,l$ which are powers of $p$,   with $k\mid l$, let 
        \[\theta_{k,l}\colon C^*(\square;\Z/k)\to C^*(\square;\Z/l)\]
        be the homomorphism induced by the inclusion $\Z/k\xrightarrow{\times l/k}\Z/l$.
        \item  Let
        \[\delta_k\colon H^*(\square;\Z/k)\to H^{*+1}(\square;\Z)\]
        denote the Bockstein homomorphism associated to the short exact sequence 
        \[0\to\Zp\xrightarrow{\times k}\Zp\to\Z/k\to 0\]
        \item Let $\pp^i$ be the conventional notation for the Steenrod reduced power operation. 
    \end{itemize}    
\end{nota}

Let $\xi_n\in C^2(\K{n},\Zp)$ such that the mod $n$ reduction of $\xi_n$ represents the fundamental class in $H^2(\K{n};\Z/n)$. Let $\eta_n\in C^3(\K{n};\Zp)$ be the cocycle such that $\delta(\xi_n) = n\eta_n$.

We define
\[\zzeta{n}{r} = \delta_n([\xi_n]^r)\in H^{2r+1}(\K{n};\Zp). \] 
For abbreviation, we write $\zeta_n = \zzeta{n}{1}$.

\begin{prop}\label{pro:Kn_integral}
    Let $n$ be as above. We have
    \begin{equation*}
        H^{k}(\K{n};\Zp)\cong
        \begin{cases}
            \Z/n,\ k= 2r + 1, 0 < r < p,\\
             0,\ k = 2r, 0 < r \leq p\\
            \Z/pn,\ k = 2p + 1,\\
           \Z/p,\ k = 2p+2.
        \end{cases}
    \end{equation*}
    For $k = 2r + 1$, $0 < r \leq p$, $H^{k}(\K{n};\Z)$ is generated by $\zzeta{n}{r}$. The group $H^{2p+2}(\K{n};\Zp)$ is generated by  $\delta_p\pp^1\rho_p([\eta_n])$. 
\end{prop}

Let $m = p^{f}$ for some $f \geq e$. 
We let 
\[\xx{n}{r} = \theta_{n,m}([\xi_n]^r),\ \yy{n}{r} = \rho_m(\zzeta{n}{r}).\]

\begin{prop}\label{pro:generators_for_HK}
    Let $m$, $n$ be as above. Then we have
    \begin{equation}
        H^k(\K{n};\Z/m)\cong
        \begin{cases}
             \Z/n,\ 2\leq k\leq 2p-1,\\
             \Z/\epsilon_p(mn^{-1})n,\ k = 2p,\\
             \Z/\epsilon_p(mn^{-1})n\oplus\Z/p,\
              k = 2p+1,\\
             \Z/n\oplus\Z/p,\ k = 2p+2.
        \end{cases}
    \end{equation}
    For $2\leq k \leq 2p$, the cyclic group $H^k(\K{n},\Z/m)$ is generated by $\xx{n}{r}$ for $k = 2r$ and by $\yy{n}{r}$ for $k = 2r+1$.
    
    The summand $\Z/\epsilon_p(mn^{-1})n$ of $H^{2p+1}(\K{n};\Z/m)$ is generated by $\yy{n}{p}$. The summand $\Z/p$ of $H^{2p+1}(\K{n};\Z/m)$ is generated by $\theta_{p,m}\pp^1\rho_p(y_n)$.
    
    The summand $\Z/p$ of $H^{2p+2}(\K{n};\Z/m)$ is generated by $\delta_p\pp^1\rho_p(y_n)$.
\end{prop}
\begin{proof}
    This follows from Section~11 of \cite{Cartan1954}.
\end{proof}

\begin{prop}\label{pro:products_in_K}
    In $H^*(\K{n};\Z/m)$, there are relations
    \begin{equation*}
        \xx{n}{i}\xx{n}{j}=\frac{m}{n}\xx{n}{i+j},\ \xx{n}{i}\yy{n}{j} = \frac{m}{n}\yy{n}{i+j}.
    \end{equation*}
\end{prop}

Let $m,n$ be as above. Let $\phi\colon\Z/n\to\Z/m$ be the injective homomorphism defined by $\phi(1) = mn^{-1}$. Then we have an induced homomorphism 
\[\phi^{\#}\colon C^*(\K{m};\Zp)\to C^*(\K{n};\Zp).\]

\begin{prop}\label{pro:Kn_to_Km}
  The cochains $\xi_n$, $\xi_m$, $\eta_n$ and $\eta_m$ may be chosen so that
 \begin{equation*}
        \phi^{\#}(\xi_m) = mn^{-1}\xi_n,\ \phi^{\#}(\eta_m) = \eta_n.
    \end{equation*}  
   
    Passing to cohomology with coefficients in $\Z/m$ or in $\Z$, we have
    \begin{equation}\label{eq:Kn_to_Km}
        \phi^*(\xx{m}{r}) = (\frac{m}{n})^{r-1}\xx{n}{r},\ \phi^*(\yy{m}{r}) = (\frac{m}{n})^{r-1}\yy{n}{r},\ \phi^*(\zzeta{m}{r}) = (\frac{m}{n})^{r-1}\zzeta{n}{r}
    \end{equation}
\end{prop}
\begin{proof}
    The suitable choices for $\xi_n$, $\xi_m$, $\eta_n$ and $\eta_m$ are easily observed in the iterated bar construction. The rest are routine computation.
\end{proof}
Let $\barxi_n,\bareta_n$ denote the mod $m$ reduction of $\xi_n$ and $\eta_n$ in $C^*(\K{n};\Z/m)$. Suppose $f = ke + e_0$ for $0\leq e_0 < e$. Let $n_0 = p^{e_0}$. Let $F^iC^*(\K{n};\Z/m) = n^iC^*(\K{n};\Z/m)$. Then $\{F^i\}_i$ defines a finite decreasing filtration on $C^*(\K{n};\Z/m)$. Let the associated spectral sequence and its differentials be denoted by $(E_r^{s,t}, d_r^{s,t})$.

The spectral sequence $(E_r^{s,t}, d_r^{s,t})$ satisfies
   \begin{equation}\label{eq:diff_mod_m_e1}
       E_1^{s,t}\cong
       \begin{cases}
          H^t(\K{n};\Z/n),\ 0\leq s\leq k-1,\\
          H^t(\K{n};\Z/n_0),\ s = k-1,\\
          0,\ \textrm{else}.
       \end{cases}
   \end{equation}
   The differentials are of the form
   \begin{equation*}
       d_r^{s,t}\colon E_r^{s,t}\to E_r^{s+r,t}.
   \end{equation*}
   To avoid ambiguity, for any $\alpha\in H^t(\K{n};\Z/n)$ or $\alpha\in H^t(\K{n};\Z/n_0)$, we let $\alpha[s]$ denote its counterpart in $E_1^{s,t}$.
\begin{lem}\label{lem:diff_mod_m}
    Let $1\leq i\leq p$.
    \begin{enumerate}[label={(\arabic*)}]
        \item\label{ite:diff_mod_m_1} For $1\leq s\leq k$, $0\leq t\leq 2p+1$, we have $E_2^{s,t} = E_{\infty}^{s,t} = 0$.
        \item\label{ite:diff_mod_m_2} For $0\leq t\leq 2p+1$, we have $E_2^{0,t} = E_{\infty}^{0,t}$ which is generated by $\barxi_n^{i-1}\bareta_n[0]$ if t =$2i+1$, and is $0$ otherwise.  
        \item\label{ite:diff_mod_m_3} For $0\leq t\leq 2p+1$, we have $E_2^{0,t} = E_{\infty}^{0,t}$ which is generated by $\barxi_n^{i}[k+1]$ if t =$2i$, and is $0$ otherwise. 
        \item\label{ite:diff_mod_m_4} The class $\rho_m(\yy{n}{i})$ is represented by a cocycle of the form
        \[\barxi_n^{i-1}\bareta_n + n\lambda.\]
        \item\label{ite:diff_mod_m_5} There is a cochain $\lambda'$ such that
        \[d(\barxi_n^i + n\lambda') = n(\barxi_n^{i-1}\bareta_n + n\lambda).\]
    \end{enumerate}
\end{lem}
\begin{proof}   
    It follows from a routine computation that we have 
    \begin{equation}\label{eq:diff_mod_m_d2}
        d_1^{s,t}(\barxi_n^i[s]) = \barxi_n^{i-1}\eta_n[s+1], \ 0\leq s\leq k,
    \end{equation}
    and the desired result follows immediately.
\end{proof}

\begin{prop}\label{pro:xyMassey}
    Let $m, n, r$ be positive integers satisfying $2\leq r \leq p-1$ and $n^r\mid m$. Then the $(r+1)$-fold Massey product
    \begin{equation}\label{eq:xyMassey}
      M_r = \langle x_n,\underbrace{y_n,\cdots,y_n}_{\textrm{$r$ tuples}}\rangle  
    \end{equation}
    is well-defined, and we have
    \begin{equation}\label{eq:xyMassey_formula}
       M_r = \{\frac{km}{r!n^r}\yy{n}{r+1}\mid k\equiv 1\pmod{n}\}. 
    \end{equation}
\end{prop}
\begin{proof}
    By induction on $r$, we show that $M_r$ is well-defined, and that any defining system for $M_r$ is of the form
    \begin{equation}\label{eq:M_defining_sys}
         \begin{matrix}
             \frac{m}{n}\barxi_n & \mmu{n}{2} &  \mmu{n}{3}&\cdots &\mmu{n}{r} & * \\
             & \bareta_n & \cdots &\cdots &\cdots &\vdots\\
             & & \ddots &\cdots & \vdots & \vdots\\
             & &        &\ddots  & \vdots & \vdots \\
             & &   &     &  \bareta_n & \vdots\\
             & &   &     &                &  \bareta_n.     
        \end{matrix}
     \end{equation}
     where 
     \[\mmu{n}{i} = \frac{k_im}{i!n^i}(\barxi_n^i + n\lambda^{(i)})\]
     where $k_i\equiv 1\pmod{n}$, and $\barxi_n^i + n\lambda^{(i)}$ represents $\xx{n}{i}$. 
     
     For $r=2$, the cochain $\frac{m}{n}\barxi_n\bareta_n$ represents the cup product $x_ny_n$. If $n^2\mid m$, then by Lemma~\ref{lem:diff_mod_m}, \ref{ite:diff_mod_m_4}, the required cochain $\mmu{n}{2}$ exists, and then the defining system \eqref{eq:M_defining_sys} exists by Lemma~\ref{lem:nor_def_sys_ext}. The indeterminacy is
     \begin{equation*}
        \begin{split}
            & x_nH^5(\K{n};\Z/m) + H^4(\K{n},\Z/m)y_n\\
           =& \{x_n\yy{n}{2} + j\xx{n}{2}y_n\} =  \{(\frac{m}{n}i +\frac{m}{n}j)\yy{n}{3} \}\\
           =& \{\frac{m}{n}i\yy{n}{3}\mid i\in\Z/m\}
        \end{split}           
     \end{equation*}
     and we have proved \eqref{eq:xyMassey_formula} for $r=2$.
         
     For the inductive argument, suppose we have $r<p-1$, $n^{r+1}\mid m$, and the defining system \eqref{eq:M_defining_sys} is given, and the Massey product 
     \[ M_r = \{\frac{km}{r!n^r}\yy{n}{r+1}\mid k\equiv 1\pmod{n}\}. \]
     Therefore, each class in $M_r$ is represented by a cocycle of the form
     \[\frac{km}{r!n^r}(\barxi_n^r\bareta_n + n\lambda),\]
     for some cochain $\lambda$ and $k\equiv1\pmod{n}$.
     
     By \eqref{eq:diff_mod_m_d2}, a cochain $\mmu{n}{r+1}$ satisfies 
     \[d(\mmu{n}{r+1}) = \frac{km}{r!n^r}(\barxi_n^i\bareta_n + n\lambda)\]
     if and only if $\mmu{n}{r+1}$ is of the form 
     \[\mmu{n}{r+1} = \frac{km}{(r+1)n^{r+1}}(\barxi_n^{r+1} + n\lambda')\]
     for some cochain $\lambda'$. By \eqref{eq:diff_mod_m_d2} again, such a cocycle $\mmu{n}{r+1}$ exists. Therefore, we have acyclic defining systems of the form 
      \begin{equation*}
         \begin{matrix}
             \frac{m}{n}\barxi_n & \mmu{n}{2} &  \mmu{n}{3}&\cdots &\mmu{n}{r} & \mmu{n}{r+1} \\
             & \bareta_n & \cdots &\cdots &\cdots &\vdots\\
             & & \ddots &\cdots & \vdots & \vdots\\
             & &        &\ddots  & \vdots & \vdots \\
             & &   &     &  \bareta_n & \vdots\\
             & &   &     &                &  \bareta_n.     
        \end{matrix}
     \end{equation*}
     By Lemma~\ref{lem:nor_def_sys_ext}, the above acyclic defining system may be augmented by a column to the right and become a defining system of the Massey product of $M_{r+1}$ and each class in $M_{r+1}$ is represented by a cocycle of the form
     \[\frac{km}{(r+1)!n^{r+1}}(\barxi_n^{r+1}\bareta_n + n\lambda'),\]
     which represents the class
     \[\frac{km}{(r+1)!n^{r+1}}\yy{n}{r+2},\ k\equiv 1\pmod{n}.\]
     The induction is therefore completed.
\end{proof}
\begin{remark}\label{rmk:first_two_mu}
A choice of $\mmu{n}{2}$ and $\mmu{n}{3}$ is given as follows:
\begin{equation*}
    \mmu{n}{2} = \frac{m}{2n^2} + \frac{m}{2n}x\cupo y,
\end{equation*}
\begin{equation*}
    \mmu{n}{3} = \frac{m}{6n^3}x^3 + \frac{m}{6n^2}[x^2\cupo y - (xy)\cupo x] + \frac{m}{6n} [(x\cupo y)\cupo y + (x\smile_2 y)y +y(x\smile_2 y)].
\end{equation*}
\end{remark}
\begin{cor}\label{cor:xyMassey}
    Let $1 \leq r \leq p-2$ and $m$, $n$ be powers of $p$ such that $n^r\mid m$. The Massey product $M_r$ is defined and contains exactly one element 
    \[\frac{m}{r!n^r}\yy{n}{r+1}\in H^{2r+3}(\K{n},\Z/m)\]
    for some $k\in(\Z/m)^{\times}$. 
    In particular, we have $M_r = 0$ if and only if $n^{r+1}\mid m$.

    If $n^{p-1}\mid m$, then $M_{p-1}$ is defined as an element in 
    \[H^{2p+1}(\K{n},\Z/m)/\frac{m}{n^{p-2}}\yy{n}{p}.\]
    and contains a class
    \[\frac{km}{(p-1)!n^{p-1}}\yy{n}{p}\]
    for some $k\in(\Z/m)^{\times}$
    If, in addition, we have $pn^{p-1}\mid m$, then $M_{p-1}$ contains exactly one class.
\end{cor}

For the rest of this section , let $m = n$ and let $p$ be an odd prime factor of $n$. For $1\leq r\leq p$ We have the generator
\[ \zzeta{n}{r}\in H^{2r+1}(K(\Z/n,2);\Zp),\]
which is of order $n$ for $1\leq r\leq p-1$ and of order $pn$ for $r = p$. The class $\zzeta{n}{p}$ defines a cohomology operation
\begin{equation}\label{eq:O_n}
    \O_r\colon H^2(-;\Z/n)\to H^{2p+1}(-;\Z).
\end{equation}
Let $\rho_{n,p}\colon H^*(-;\Z/n)\to H^*(-;\F_p)$ be the mod $p$ reduction and recall the Bockstein homomorphism $\delta_p$ as in Notation~\ref{not:operations}. 

\begin{prop}\label{pro:O_n}
    For a topological space $X$ and $w\in H^2(X;\Z/n)$, we have
    \[\delta_p(\rho_{n,p}(w)^p) = n\O_p(w).\]    
\end{prop}



\section{The cohomology of $BPU(m)$}\label{sec:cohomology}
This section is a brief review of the cohomology of $BPU(m)$, for $p\mid m$. We use the standard notations $\pp^i$ for the Steenrod operations and $\delta_p\colon H^*(\Box;\Z/p)\to H^{*+1}(\Box;\Z)$ for the Bockstein homomorphism.

An elementary calculation shows 
\begin{equation}\label{eq:H3}
    H^k(BPU(m);\Z)\cong
    \begin{cases}
        0,\ k = 1,2\\
        \Z/m,\ k = 3.
    \end{cases}
\end{equation}
There is a canonical generator of $H^3(BPU(m);\Z)$ which we denote by $\mathfrak{b}$.

We summarize what we need about the cohomology of $BPU(m)$ as follows.
\begin{prop}[\cite{gu2021p}]\label{pro:BPU}
    Let $p$ be an odd prime and $m$ be a power of $p$. Then we have 
    \begin{equation}\label{eq:BPU_cases}
        H^k(BPU(m);\Zp)\cong
        \begin{cases}
            0,\ k = 1,2,\\
             (\Z/m)_{(p)},\ k = 3,\\
            0,\ 4 \leq k < 2p+2, k\textrm{ is odd}, \\
            F_k,\ 4 \leq k < 2p+2, k\textrm{ is even},
            F_k\oplus\Z/p,\ k = 2p +2,
        \end{cases}
    \end{equation}
   where $F_k$ is a finitely generated free $\Zp$-module. The summand $\Z/p$ of $H^{2p+2}(BPU(m);\Zp)$ is generated by $\delta_p\pp^1\rho_p(\mathfrak{b})$.
\end{prop}

Let $\uinf{m}\in H^{2p+1}(BPU(m);\Z/m)$ be the image of $\pp^1\rho_p(\mathfrak{b})$ via the homomorphism
\[H^{2p+1}(BPU(m);\Z/p)\to H^{2p+1}(BPU(m);\Z/m)\]
induced by the inclusion $\Z/p\to\Z/m$. Let
\[\vinf{m} = \rho_m(\delta_p\pp^1\rho_p(\mathfrak{b}))\in H^{2p+1}(BPU(m);\Z/m).\]
The seemingly peculiar choice of notations is justified in Section~\ref{sec:existence}. A routine computation yields the following
\begin{cor}\label{cor:BPU}
    We have $H^{2p+1}(BPU(m);\Z/m)\cong\Z/p$, generated by $\uinf{m}$.
\end{cor}

\section{The existence of certain principal $PU(m)$-bundles}\label{sec:existence}
Let $m$ be as in Section~\ref{sec:Kn2}. By Bott periodicity, the homotopy groups of $BPU(m)$ satisfies
\begin{equation}\label{eq:htpy_grp}
    \pi_k(BPU(m))\cong
    \begin{cases}
        \Z/m,\ k=2,\\
        \Z, \ k\textrm{ even},\ k < 2m + 1,\\
        0, \ k\textrm{ odd},\ k < 2m + 1.    
    \end{cases}
\end{equation}

Let $X_k$ be the $k$th level of the Postnikov tower of $BPU(m)$. By \eqref{eq:htpy_grp}, we have $X_{2r} = X_{2r+1}$ for $r < p + 1$ and $X_2 = \K{m}$. For $r\geq 1$, let $\htf{2r}$ denote the map $X_{2r}\to X_{2(r-1)}$, and let $f_{2r}$ denote the composition
\[X_{2r}\xrightarrow{\htf{2r}} X_{2(r-1)}\cdots \to X_2 = \K{m}.\]
In particular $f_{2} = \htf{2} = \opn{id}_{X_2}$. 

By \eqref{eq:htpy_grp}, we have a homotopy fiber sequence
\begin{equation}\label{eq:Xfiber_seq}
    X_{2r+2}\xrightarrow{\htf{2r+2}} X_{2r}\xrightarrow{\kappa_{2r+3}}K(\Z,2r+3)
\end{equation}
where $\kappa_{2r+3}$ is the $(2r+3)$th k-invariants of the Postikov tower of $BPU(m)$.

We proceed to define a collection of cocycles and cohomology classes for $X_{2r}$ based on the notations introduced in Proposition~\ref{pro:Kn_integral} and Proposition~\ref{pro:generators_for_HK}.

Let $\xxi{m}{r} = f_{2r}^{\#}(\xxi{m}{1})$, $\eeta{m}{r} = f_{2r}^{\#}(\eeta{m}{1})$, and let 
\[\xxx{m}{r}{k} = [\xxi{m}{r}^k]\in H^{2k}(X_{2r};\Z/m), \ \yyy{m}{r}{k} = [\xxi{m}{r}^{k-1}\eeta{m}{r}]\in H^{2k+1}(X_{2r};\Z/m).\]
Let $\ssig{m,r}{k} = f_{2r}^{\#}(\ssig{m}{k})$, $\ttau{m,r}{k} = f_{2r}^{\#}(\ttau{m}{r})$, and let 
\[u_{m,r}^{(k)} = [mp^{-1}\sigma_{m,r}^{(k)}]\in H^{2p^k+1}(X_{2r};\Z/m),\ v_{m,r}^{(k)} = [\ttau{m,r}{k}]\in H^{2p^k+2}(X_{2r};\Z/m).\]
Finally, let 
\[\zzzeta{m}{r}{k} = [\xxi{m}{r}^{k-1}\eeta{m}{r}]\in H^{2k+1}(X_{2r};\Zp).\]
For abbreviation, let $x_{m,r} = \xxx{m}{r}{1}$, $y_{m,r} = \yyy{m}{r}{1}$, $\zeta_{m,r} = \zzzeta{m}{r}{1}$, $u_{m,r} = u_{m,r}^{(1)}$ and $v_{m,r} = v_{m,r}^{(1)}$.

We have a variation of \eqref{eq:Xfiber_seq}:
    \begin{equation}\label{eq:inductive_fiber_seq}
        K(\Z, 2r)\to X_{2r}\xrightarrow{\htf{2r}} X_{2r-2}.
    \end{equation}
Let $(\wE_*^{*,*}, \wdd_*^{*,*})$ be its associated cohomological Serre spectral sequences with coefficients in $\Z_{(p)}$.
Let $(E_*^{*,*}, d_*^{*,*})$ be the one with coefficients in $\Z/m$.

Let $e_{2r}\colon BPU(m)\to X_{2r}$ be the canonical map.

For a $\Zp$-module $M$, let $T(M)$ denote the submodule of torsion elements of $M$. 
\begin{lem}\label{lem:rational}
    For $r\geq 1$, $k$ odd, we have
    \[T(H^k(X_{2r};\Zp)) = 0.\]
\end{lem}
\begin{proof}
    By construction, we have
    \begin{equation*}
        \pi_i(X_{2r})\cong
        \begin{cases}
            \Z,\ i \leq 2r,\ i\textrm{ is even},\\
            0,\ \textrm{otherwise.}
        \end{cases}
    \end{equation*}
    Therefore, $H^*(X_{2r};\mathbb{Q})$ is a polynomial ring concentrated in even dimensions, and we conclude.
\end{proof}

Recall that $\rho_m\colon H^*(-;\Z)\to H^*(-;\Z/m)$ denote the mod $m$ reduction homomorphism.

\begin{lem}\label{lem:2p+1_X}
    Let $r\geq 1$. We have $u_{m,r}\in H^{2p+1}(X_{2r};\Z/m)$ such that 
    \[\delta_m(u_{m,r}) = \delta_p\pp^1\rho_p(y_{m,r}) = v_{m,r}\neq 0.\]
    Moreover, we have $f_{2r}^*(u_{m}) = u_{m,r}$ and $e^*_{2r}(u_{m,r}) = \uinf{m}$.
\end{lem}

\begin{lem}\label{lem:stable}
    Let $r\geq 1$.
    \begin{enumerate}
        \item\label{ite:stable_X}
        The induced homomorphism $\htf{2r}^*\colon H^k(X_{2r-2};\Zp)\to H^k(X_{2r};\Zp)$ is an isomorphism for $k < 2r$ or $k=2r+2$, a split injection with cokernel $\Zp$ for $k = 2r$, and a surjection for $k = 2r+1$. 
        \item\label{ite:stable_bpu}   The induced homomorphism $e_k^*\colon H^k(X_{2r};\Zp)\to H^k(BPU(m);\Zp)$ is an isomorphism for $k < 2r + 2$, and a split injection with cokernel $\Zp$ for $k = 2r+2$.
    \end{enumerate}
\end{lem}
\begin{proof}
    We prove (\ref{ite:stable_X}). In the spectral sequence $(\wE,\wdd)$, for $s< 2r$, for degree reasons we have 
    \[H^s(X_{2r-2};\Zp)\cong E_{2}^{s,0}=E_{\infty }^{s,0} \cong H^s(X_{2r};\Zp).\]
    Therefore, we have
    \begin{equation}\label{eq:stable_1}
         H^s(X_{2r-2};\Zp)\xrightarrow[\cong]{\htf{2r}^*}H^s(X_{2r};\Zp),\ s < 2r.
    \end{equation}
    We have a short exact sequence
    \begin{equation}\label{eq:stable_2}
        0\to \wE_{\infty}^{2r,0}\to H^{2r}(X_{2r};\Zp)\to \wE_{\infty}^{0,2r}\to 0.
    \end{equation}
    By Lemma~\ref{lem:rational}, the differential $\wdd_{2r+1}^{0,2r}$ has domain $E_2^{0,2r}\cong\Zp$ and image a torsion submodule. Therefore, we have $E_{\infty}^{0,2r}\cong\Zp$. On the other hand, for degree reasons we have $\wE_2^{2r,0} = \wE_{\infty}^{2r,0}$. Therefore, the exact sequence \eqref{eq:stable_2} becomes
    \begin{equation}\label{eq:stable_3}
        0\to H^{2r}(X_{2r-2};\Zp)\xrightarrow{\htf{2r}^*} H^{2r}(X_{2r};\Zp)\to \Zp\to 0.         
    \end{equation}
    Again, for degree reasons, we have 
    \begin{equation}\label{eq:stable_4}
       H^{2r+1}(X_{2r};\Zp) \cong\wE_{\infty}^{2r+1,0} = \wE_2^{2r+1,0}/\opn{Im}\wdd_{2r+1}^{0,2r},        
    \end{equation}
    and
    \begin{equation}\label{eq:stable_5}
        E_{\infty}^{2r+2,0} = E_2^{2r+2,0} = H^{2r+2}(X_{2r-2};\Zp)\xrightarrow[\cong]{\htf{2r}^*}H^{2r+2}(X_{2r};\Zp).        
    \end{equation}
    Then (\ref{ite:stable_X}) follows from \eqref{eq:stable_1}, \eqref{eq:stable_3},  \eqref{eq:stable_4} and \eqref{eq:stable_5}.
    
    We prove (\ref{ite:stable_bpu}). By replacing $r$ with $r+1, r+2, \cdots$ in \eqref{eq:stable_1}, we obtain     
    \begin{equation*}
        H^s(X_{2r-2};\Zp)\xrightarrow[\cong]{e^*_{2r}}H^s(BPU(m);\Zp),\ s< 2r.
    \end{equation*}
    Replacing $r$ with $r+1$, we have
    \begin{equation}\label{eq:stable_6}
        H^s(X_{2r};\Zp)\xrightarrow[\cong]{e^*_{2r+2}}H^s(BPU(m);\Zp),\ s< 2r+2.
    \end{equation}
    Replace $r$ with $r+1$ in \eqref{eq:stable_3} and \eqref{eq:stable_6}, we obtain (\ref{ite:stable_bpu}). 
\end{proof}

\begin{lem}\label{lem:cyclic} 
    Let $1\leq r\leq p-1$.
    \begin{enumerate}
    \item\label{ite:cyclic1} The $\Zp$-module $H^{2r+3}(X_{2r};\Zp)$ is finite cyclic and generated by $\kappa_{2r+3}$, the k-invariant of  $BPU(m)$. 
    
    \item\label{ite:cyclic2}There is an exact sequence
        \begin{equation}\label{eq:SES_int}
            0\to \wE_{\infty}^{2r+3,0}\to H^{2r+3}(X_{2r};\Zp)\to\wE_{\infty}^{3,2r}\to 0,
        \end{equation}
        and we have $\wE_{\infty}^{2r+3,0} = \wE_2^{2r+3,0}$.
    \end{enumerate}    
\end{lem}
\begin{proof}
    We prove (\ref{ite:cyclic1}). By construction, the kernel of 
    \[\htf{2r+2}^*\colon H^{2r+3}(X_{2r};\Zp)\to H^{2r+3}(X_{2r+2};\Zp)\]
    is generated by $\kappa_{2r+3}$. By Lemma~\ref{lem:stable}, (\ref{ite:2r+4_stable}), we have 
    \[H^{2r+3}(X_{2r+2};\Zp)\xrightarrow[\cong]{e_{2r+2}^*} H^{2r+3}(BPU(m);\Zp).\]
    By Proposition~\ref{pro:BPU}, we have $H^{2r+3}(BPU(m);\Zp) = 0$.

    Combining the arguments above, we conclude that $H^{2r+3}(X_{2r+2};\Zp)$ is generated by $\kappa_{2p+3}$, which is a finite cyclic group, by Lemma~\ref{lem:rational}.
    
    We prove (\ref{ite:cyclic2}). The exact sequence \eqref{eq:SES_int} follows from the observation that the only nontrivial entries on the $\wE_2$-page of total degree $2r+3$ are $\wE_2^{3,2r}$ and $\wE_2^{2r+3,0}$.

    By comparing degrees, it follows that $\wE_*^{2r+3,0}$ supports and receives no nontrivial differential. Therefore, we have $\wE_2^{2r+3,0} = \wE_{\infty}^{2r+3,0}$. 
\end{proof}

\begin{lem}\label{lem:cyclic_m}
        For $1 \leq r \leq p-1$, we have a short exact sequence
         \begin{equation}\label{eq:SES_modm}
            0\to E_{\infty}^{2r+3,0}\to H^{2r+3}(X_{2r};\Z/m)\to E_{\infty}^{3,2r}\to 0,
        \end{equation}
        and we have $E_{\infty}^{2r+3,0} = E_2^{2r+3,0}$.
\end{lem}
\begin{proof}
    The only nontrivial entries on the $E_2$-page of total degree $2r+3$ are $E_2^{3,2r}$ and $E_2^{2r+2,0}$, and the short exact sequence follows. The equation $E_{\infty}^{2r+3,0} = E_2^{2r+3,0}$ follows for obvious degree reasons.
\end{proof}

\begin{lem}\label{lem:2r+4}
    \begin{enumerate}
        \item\label{ite:2r+4_r}         
        For $r < p-1$, we have 
        \[T(H^{2r+4}(X_{2r-2};\Zp)) = T(\wE_2^{2r+4,0}) =  \opn{Im}\wdd_{2r+1}^{3,2r}.\]
        \item\label{ite:2r+4_p-1} We have 
        \[T(H^{2p+2}(X_{2p-4};\Zp)/\opn{Im}\wdd_{2p-1}^{3,2p-2})\cong \Z/p\]
        which is generated by the coset represented by $v_{m,p-2}$.
        \item\label{ite:2r+4_stable} For $r\leq p-1$, we have induced isomorphisms
        \[T(H^{2r+4}(X_{2r-2};\Zp)/\opn{Im}\wdd_{2r+1}^{3,2r})\xrightarrow[\cong]{\htf{2r}^*}T(H^{2r+4}(X_{2r};\Zp))\xrightarrow[\cong]{e_{2r}^*}T(H^{2r+4}(BPU(m);\Zp)).\]
    \end{enumerate}
\end{lem}
\begin{proof}
    By Lemma~\ref{lem:rational}, we have    $\opn{Im}\wdd_{2r+1}^{3,2r}\subset T(\wE_2^{2r+4,0})$. Therefore, we have 
    \[T(H^{2r+4}(X_{2r-2};\Zp)/\opn{Im}\wdd_{2r+1}^{3,2r}) = T(H^{2r+4}(X_{2r-2};\Zp))/\opn{Im}\wdd_{2r+1}^{3,2r}.\]

    Let $0\leq r\leq p-1$. 
    In the spectral sequence $(\wE_{*}^{*,*}, \wdd_*^{*,*})$, the only non-torsion-free entry of total degree $2r+4$ is $\wE_*^{2r+4,0}$. For degree reasons, we have 
    \begin{equation}\label{eq:2r+4_1}
        \begin{split}
            & T(H^{2r+4}(X_{2r};\Zp)) \cong T(\wE_{\infty}^{2r+4,0}) = T(\wE_{2r+2}^{2r+4,0})\\
            = & T(\wE_2^{2r+4,0}/\opn{Im}\wdd_{2r+1}^{3,2r})\cong T(H^{2r+4}(X_{2r-2};\Zp)/\opn{Im}\wdd_{2r+1}^{3,2r}).
        \end{split}
    \end{equation}

    By Lemma~\ref{lem:stable}, (\ref{ite:stable_X}), (\ref{ite:stable_bpu}), we have
    \begin{equation}\label{eq:2r+4_2}
        T(H^{2r+4}(X_{2r};\Zp))\xrightarrow[\cong]{f_{2r+2}^*} T(H^{2r+4}(X_{2r+2};\Zp))\xrightarrow[\cong]{e_{2r+2}^*} T(H^{2r+4}(BPU(m);\Zp)).
    \end{equation}
    Combining \eqref{eq:2r+4_1} and \eqref{eq:2r+4_2}, we have
    \begin{equation}\label{eq:2r+4_3}
        T(H^{2r+4}(X_{2r-2};\Zp)/\opn{Im}\wdd_{2r+1}^{3,2r})\xrightarrow[\cong]{e_{2r-2}^*}T(H^{2r+4}(BPU(m);\Zp))\cong
        \begin{cases}
            0,\ r< p-1,\\
            \Z/p,\ r = p-1.
        \end{cases}
    \end{equation}
    This proves (\ref{ite:2r+4_r}) and (\ref{ite:2r+4_p-1}). For (\ref{ite:2r+4_stable}), notice that the isomorphism between both ends of \eqref{eq:2r+4_1} is induced by $\htf{2r+2}^*$, and that $e_{2r}^* = e_{2r+2}^*\htf{2r+2}^*$ is an isomorphism between the torsion submodules, by \eqref{eq:2r+4_2}.  
\end{proof}

For the next two propositions, recall the convention $\bar{\alpha} = (-1)^{\opn{deg}\alpha}\alpha$ where $\alpha$ is a homogeneous element in a cochain complex.

\begin{prop}\label{pro:Massey_in_X}
   For $1\leq r\leq p-2$, or $r = p-1$, $p^2\mid m$, the $(r+1)$-fold Massey product
   \[L_{m,r} = \langle x_{m,r},\underbrace{y_{m,r},\cdots, y_{m,r}}_{r\textrm{ tuples}}\rangle\]
   is well-defined, and has a defining system of the form
    \begin{equation*}
         \begin{matrix}
             \xxi{m}{r} & \om{r}{1}{1} & \om{r}{1}{2}&\cdots & \om{r}{1}{r-1} & * \\
             & \eeta{m}{r} & \om{r}{2}{2} &\cdots &\om{r}{2}{r-1} & \om{r}{2}{r}\\
             & & \ddots &\cdots & \vdots & \vdots\\
             & &   &     &  \eeta{m}{r} & \om{r}{r}{r}\\
             & &   &     &                &  \eeta{m}{r}     
        \end{matrix}
     \end{equation*}  
    of which the second to the last rows form an $r$-fold normal acyclic defining system for the $r$-fold Massey product 
     \[\langle\underbrace{y_{m,r},\cdots, y_{m,r}}_{r\textrm{ tuples}}\rangle,\] 
     and there is an integer $q_r$, coprime to $p$, such that    
     \begin{enumerate}
         \item\label{ite:Massey_in_X_r}  for $1\leq r \leq p-2$, we have 
             \begin{equation*}
                \rho_m(\kappa_{2r+3}) =  q_r[\bxxi{m}{r}\om{r}{2}{r} + \sum_{i=1}^{r-1}\bom{r}{1}{i}\om{r}{i+1}{r} + \bom{r}{1}{r-1}\eeta{m}{r}],
            \end{equation*}
            which generates $H^{2r+3}(X_{2r};\Z/m)$, and 
        \item\label{ite:Massey_in_X_r_p-1}  for $r = p-1$ and $p^2\mid m$, we have
        \begin{equation*}
            \rho_m(\kappa_{2p+1}) \equiv  q_{p-1}[\bxxi{m}{p-1}\om{p-1}{2}{p-1} + \sum_{i=1}^{p-2}\bom{p-1}{1}{i}\om{p-1}{i+1}{p-1} + \bom{p-1}{1}{p-2}\eeta{m}{p-1}]\pmod{u_{m,p-1}},
        \end{equation*}
        which generates $H^{2r+3}(X_{2r};\Z/m)/(u_{m,p-1})$.
    \end{enumerate}  
\end{prop}
\begin{proof}    
   We proceed by induction in $r$. In the case $r = 1$, we have $X_{2r} = X_2 = \K{m}$, and the desired result follows from \eqref{eq:iota} and the relation $\rho(\zzeta{m}{2}) = x_my_m$.  This is also shown by Antieau-Williams (\cite{antieau2014topological}, Proposition~4.4).

   For $r > 1$, by the induction hypothesis,the $r$-fold Masssey product
   \begin{equation}\label{eq:inductive_hyp}
      L_{m,r-1} = \langle x_{m,r-1},\underbrace{y_{m,r-1},\cdots, y_{m,r-1}}_{r-1\textrm{ tuples}}\rangle,
   \end{equation}
   is defined by a defining system of the form
   \begin{equation}\label{eq:def_sys_r-1}
        \begin{matrix}
             \xxi{m}{r-1} & \om{r-1}{1}{1} & \om{r-1}{1}{2}&\cdots & \om{r-1}{1}{r-2} & * \\
             & \eeta{m}{r-1} & \om{r-1}{2}{2} &\cdots &\om{r-1}{2}{r-2} & \om{r-1}{2}{r-1}\\
             & & \ddots &\cdots & \vdots & \vdots\\
             & &   &     &  \eeta{m}{r-1} & \om{r-1}{r-1}{r-1}\\
             & &   &     &                &  \eeta{m}{r-1}.     
        \end{matrix}
    \end{equation}
    of which the second to last rows form an $(r-1)$-fold normal acyclic defining system, and satisfies
    \begin{equation}\label{eq:rho_kappa_r-1}
         \rho_m(\kappa_{2r+1}) = q_{r-1}[\bxxi{m}{r-1}\om{r-1}{2}{r-1} + \sum_{i=1}^{r-3}\bom{r-1}{1}{i}\om{r-1}{i+2}{r-1} + \bom{r-1}{1}{r-2}\eeta{m}{r-1}]
    \end{equation}
    for some $q_{r-1}$ coprime to $p$. 
    
    By Lemma~\ref{lem:nor_def_sys_ext}, there are $\om{r-1}{j}{r}$, $2\leq j\leq r$ such that there is an $r$-fold normal acyclic defining system of the form
    \begin{equation}\label{eq:nor_ac_def_y}
        \begin{matrix}
            & \eeta{m}{r-1} & \om{r-1}{2}{2} &\cdots &\om{r-1}{2}{r-1} & \om{r-1}{2}{r}\\
             & & \ddots &\cdots & \vdots & \vdots\\
             & &   &     &  \eeta{m}{r-1} & \om{r-1}{r}{r}\\
             & &   &     &                &  \eeta{m}{r-1}.
        \end{matrix}
    \end{equation}
    
    We apply the chain map induced by $X_{2r}\xrightarrow{\htf{2r}} X_{2(r-1)}$ to \eqref{eq:def_sys_r-1} and obtain a defining system for $\htf{2r}^*(L_{m,r-1})$ of the form 
    \begin{equation}\label{eq:def_sys_r}
        \begin{matrix}
             \xxi{m}{r} & \om{r}{1}{1} & \om{r}{1}{2}&\cdots & \om{r}{1}{r-2} & * \\
             & \eeta{m}{r} & \om{r}{2}{2} &\cdots &\om{r}{2}{r-2} & \om{r}{2}{r-1}\\
             & & \ddots &\cdots & \vdots & \vdots\\
             & &   &     &  \eeta{m}{r} & \om{r}{r-1}{r-1}\\
             & &   &     &                &  \eeta{m}{r}.     
        \end{matrix}
     \end{equation}  
     We apply $\htf{2r}^*$ to \eqref{eq:nor_ac_def_y} and obtain an $r$-fold normal acyclic defining system for the $r$-fold Massey product
     \[\langle\underbrace{y_{m,r},\cdots, y_{m,r}}_{r\textrm{ tuples}}\rangle\]
    as follows:  \begin{equation}\label{eq:nor_ac_def_y_r}
        \begin{matrix}
            & \eeta{m}{r} & \om{r}{2}{2} &\cdots &\om{r}{2}{r-1} & \om{r}{2}{r}\\
             & & \ddots &\cdots & \vdots & \vdots\\
             & &   &     &  \eeta{m}{r} & \om{r}{r}{r}\\
             & &   &     &                &  \eeta{m}{r}.
        \end{matrix}
    \end{equation}
    Let $\iota_{2r}$ be the canonical generator of $H^{2r}(K(\Z, 2r);\Z)\cong\wE_2^{0,2r}$. Then $\iota_{2r}$ is transgressive and we have
    \begin{equation}\label{eq:iota}
        \wdd_{2r+1}^{0,2r}(\iota_{2r}) = q\kappa_{2r+1}      
    \end{equation}
    for some $q$ coprime to $p$.
    Taking mod-$m$ reduction, we have
    \begin{equation}\label{eq:iota_mod_m}
        d_{2r+1}^{0,2r}(\rho_m(\iota_{2r})) = \rho_m(q\kappa_{2r+1}).
    \end{equation}
     
    Let $d$ denote the coboundary operator in $C^*(X_{2r};\Z/m)$. By \eqref{eq:iota_mod_m}, we have $\om{r}{1}{r-1}\in C^{2r}(X_{2r};\Z/m)$ that survives to $E_2$-page, and satisfies 
     \begin{equation}\label{eq:omega_iota}
        [\om{r}{1}{r-1}] = q^{-1}q_{r-1}^{-1}\iota_{2r}\in E_2^{0,2r}=E_{2r+1}^{0,2r},   
     \end{equation}
    and    \begin{equation}\label{eq:d_omega}
         d(\om{r}{1}{r-1})  = \bxxi{m}{r}\om{r}{2}{r-1} + \sum_{i=1}^{r-3}\bom{r}{1}{i}\om{r}{i+2}{r-1} + \bom{r}{1}{r-2}\eeta{m}{r}.
     \end{equation}

    By \eqref{eq:def_sys_r}, \eqref{eq:nor_ac_def_y_r} and \eqref{eq:d_omega}, we have a defining system of the $(r+1)$-fold Massey product $L_{m,r}$ of the form
      \begin{equation}\label{eq:def_sys_L_{m,r}}
         \begin{matrix}
             \xxi{m}{r} & \om{r}{1}{1} & \om{r}{1}{2}&\cdots & \om{r}{1}{r-1} & * \\
             & \eeta{m}{r} & \om{r}{2}{2} &\cdots &\om{r}{2}{r-1} & \om{r}{2}{r}\\
             & & \ddots &\cdots & \vdots & \vdots\\
             & &   &     &  \eeta{m}{r} & \om{r}{r}{r}\\
             & &   &     &                &  \eeta{m}{r}.     
        \end{matrix}
     \end{equation}
     of which the second to the last rows form an $r$-fold acyclic defining system, and we have
     \begin{equation}\label{eq:L_Massey}
         [\bxxi{m}{r}\om{r}{2}{r} +\sum_{i=1}^{r-2}\bom{r}{1}{i}\om{r}{i+2}{r} + \bom{r}{1}{r-1}\eeta{m}{r}]\in L_{m,r}.
     \end{equation}
    Furthermore, we have
    \begin{equation}\label{eq:image_h_2r}
        \bxxi{m}{r}\om{r}{2}{r} +\sum_{i=1}^{r-2}\bom{r}{1}{i}\om{r}{i+2}{r} = \htf{2r}^{\#}(\bxxi{m}{r-1}\om{r-1}{2}{r} +\sum_{i=1}^{r-2}\bom{r-1}{1}{i}\om{r-1}{i+2}{r}). 
    \end{equation}
    By \eqref{eq:L_Massey} and \eqref{eq:image_h_2r}, the cochain complex 
    \begin{equation}\label{eq:cocycle}
        \htf{2r}^{\#}(\bxxi{m}{r-1}\om{r-1}{2}{r} +\sum_{i=1}^{r-2}\bom{r-1}{1}{i}\om{r-1}{i+2}{r}) + \bom{r}{1}{r-1}\eeta{m}{r}\in Z^{2r+3}(X_{2r};\Z/m)
    \end{equation}
    is a cocycle. 
   Furthermore, we have $\bom{r}{1}{r-1}\eeta{m}{r}$ generates $E_2^{3,2r}$, and by Lemma~\ref{lem:cyclic_m}, we have
   \begin{equation}\label{eq:E2=Einfty_modm}
       E_2^{3,2r} = E_{\infty}^{3,2r}\textrm{ generated by an element of }L_{m,r}.
   \end{equation}
    In particular, this implies $d_{2r-1}^{3,2r} = 0$, which implies 
    \begin{equation}\label{eq:wdd_2r_mod_m}
        \opn{Im}\wdd_{2r+1}^{3,2r}\subset m\wE_2^{2r+4,0} = mH^{2r+4}(X_{2r-2};\Zp).
    \end{equation}

    We prove (\ref{ite:Massey_in_X_r}). Let $r < p-1$.
    Combining \eqref{eq:wdd_2r_mod_m} and Lemma~\ref{lem:2r+4}, (\ref{ite:2r+4_r}), which states $T(\wE_2^{2r+4,0}) = \opn{Im}\wdd_{2r+1}^{3,2r}$, we have
    \[T(\wE_2^{2r+4,0})\subset mT(\wE_2^{2r+4,0}),\]
    and therefore
    \begin{equation}\label{eq:T2r-2}
      \opn{Im}\wdd_{2r+1}^{3,2r} = T(\wE_{2}^{2r+4,0}) \cong T(H^{2r+4}(X_{2r-2};\Zp)) = 0,
   \end{equation}  
   and  
   \begin{equation}\label{eq:E_2=Einfty}
        \wE_{\infty}^{3,2r} = \wE_2^{3,2r}.
   \end{equation}
   By Lemma~\ref{lem:2r+4}, (\ref{ite:2r+4_stable}), we have 
   \begin{equation}\label{eq:T2r+4}
       T(H^{2r+4}(X_{2r-2};\Zp)/\opn{Im}\wdd_{2r+1}^{3,2r})\xrightarrow[\cong]{\htf{2r}^*}T(H^{2r+4}(X_{2r};\Zp)).
   \end{equation}
   By \eqref{eq:T2r-2} and \eqref{eq:T2r+4}, we have 
   \begin{equation}\label{H2r+4_torsion_free}
        T(H^{2r+4}(X_{2r};\Zp)) = 0.     
   \end{equation}
   Therefore, we have the following surjection:
   \begin{equation}\label{eq:rho_surj}
        \rho_m\colon H^{2r+3}(X_{2r};\Zp)\twoheadrightarrow H^{2r+3}(X_{2r};\Z/m).       
   \end{equation}
   By \eqref{eq:image_h_2r}, the element 
   \[[\bxxi{m}{r}\om{r}{2}{r} +\sum_{i=1}^{r-2}\bom{r}{1}{i}\om{r}{i+1}{r} + \bom{r}{1}{r-1}\eeta{m}{r}]\in H^{2r+3}(X_{2r};\Z/m)\]
   represents $[\bom{r}{1}{r-1}\eeta{m}{r}]\in E_{\infty}^{3,2r}$. By \eqref{eq:omega_iota} and \eqref{eq:E_2=Einfty}, we have
   \begin{equation}\label{eq:identity_in_E3}
       [\om{r}{1}{r-1}\eeta{m}{r}] = q^{-1}q_{r-1}^{-1}\iota_{2r}y_{m,r}\in E_2^{3,2r} = E_{\infty}^{3,2r}
   \end{equation}
   which generates $E_{\infty}^{3,2r}\cong H^{2r+3}(X_{2r};\Z/m)$. By \eqref{eq:rho_surj} and Lemma~\ref{lem:cyclic}, (\ref{ite:cyclic1}),  $\rho_m(\kappa_{2r+3})$ is another generator of $H^{2r+3}(X_{2r};\Z/m)$, and we have an integer $q_r$, coprime to $p$ such that 
    \[\rho_m(\kappa_{2r+3}) = q_r[\bxxi{m}{r}\om{r}{2}{r} +\sum_{i=1}^{r-2}\bom{r}{1}{i}\om{r}{i+1}{r} + \bom{r}{1}{r-1}\eeta{m}{r}],\]
    and the proof of (\ref{ite:Massey_in_X_r}) is concluded. 

    We prove (\ref{ite:Massey_in_X_r_p-1}). Suppose we have $p^2\mid m$. By \eqref{eq:wdd_2r_mod_m} and Lemma~\ref{lem:2r+4}, (\ref{ite:2r+4_p-1}), we have a surjection
    \[\Z/p\cong T(\wE_2^{2p+2,0})/\opn{Im}\wdd_{2r-1}^{3,2p-2}\twoheadrightarrow T(\wE_2^{2p+2,0})/(m).\]
    So we have either $T(\wE_2^{2p+2,0})/(m) \cong \Z/p$ or $T(\wE_2^{2p+2,0})/(m) = 0$. But the latter implies $T(\wE_2^{2p+2,0}) = 0$, which contradicts Lemma~\ref{lem:2r+4}, (\ref{ite:2r+4_p-1}). Therefore, we have $T(\wE_2^{2p+2,0})/(m) \cong \Z/p$. Since $p^2\mid m$, this implies that
    \begin{equation}\label{eq:T_Zmodp}
        T(\wE_2^{2p+2,0})\cong \Z/p, 
    \end{equation}
    which, by Lemma~\ref{lem:2r+4}, (\ref{ite:2r+4_p-1}), is generated by $\nu_{m,p-1}$. Combining \eqref{eq:T_Zmodp} and $p^2\mid m$, we have $mT(\wE_2^{2p+2,0}) = 0$, and therefore 
    \begin{equation}\label{eq:wd=0_p-1}
        \wdd_{2p-1}^{3,2p-2} = 0,\ \wE_2^{3,2p-2} = \wE_{\infty}^{3,2p-2}.        
    \end{equation}
    By Lemma~\ref{lem:cyclic_m} and \eqref{eq:wd=0_p-1}, we have a short exact sequence 
    \begin{equation}\label{eq:ses2p+2_mod_m}
        0\to E_2^{2p+1,0}\to H^{2p+1}(X_{2p-2};\Zp)\to E_2^{3,2p-2}\to 0.        
    \end{equation}
    By Lemma~\ref{lem:cyclic}, (\ref{ite:cyclic1}), $\wE_{\infty}^{3,2p-2} = \wE_2^{3,2p-2}$ is generated by $\kappa_{2p+3}$. Therefore, $\rho_m(\kappa_{2p+3})$ generates $E_2^{3,2p-1}$, and by \eqref{eq:ses2p+2_mod_m}, we have $q_{p-1}$ invertible such that 
    \begin{equation}\label{eq:modEinfty2p+1}
        \rho_m(\kappa_{2p+1}) \equiv  q_{p-1}[\bxxi{m}{r}\om{r}{2}{p-1} + \sum_{i=1}^{p-2}\bom{p-1}{1}{i}\om{p-1}{i+1}{p-1} + \bom{p-1}{1}{p-2}\eeta{m}{p-1}]\pmod{E_{\infty}^{2p+1,0}}.
    \end{equation}
    By \eqref{eq:wd=0_p-1} and Lemma~\ref{lem:stable}, (\ref{ite:stable_bpu}), we have 
    \begin{equation}\label{eq:T_p-1}
         T(H^{2p+2}(X_{2p-4};\Zp))\xrightarrow[\cong]{\htf{2p-2}^*} T(H^{2p+2}(X_{2p-2};\Zp))\xrightarrow[\cong]{e_{2p-2}^*} T(H^{2p+2}(BPU(m);\Zp)) ,
    \end{equation}
    and by Proposition~\ref{pro:BPU}, we have 
    \begin{equation}\label{eq:T2p+2}
        T(H^{2p+2}(X_{2p-2};\Zp))\cong\Z/p,\ \textrm{generated by } \nu_{m,p-1}.        
    \end{equation}
    Therefore, we have 
    \begin{equation}\label{eq:2p+1_mod_m}
        \begin{split}
            & H^{2p+1}(X_{2p-2};\Z/m) = \rho_m(H^{2p+1}(X_{2p-2};\Zp))\oplus\langle u_{m,p-1}\rangle\\
        = &\langle\rho_m(\kappa_{2p+1})\rangle\oplus\langle u_{m,p-1}\rangle.
        \end{split}            
    \end{equation}
    By \eqref{eq:wd=0_p-1} and Lemma~\ref{lem:cyclic}, (\ref{ite:cyclic2}), we have the split long exact sequence
    \begin{equation}\label{eq:SES_2p+3}
        0\to E_{\infty}^{2p+1,0}\to H^{2p+1}(X_{2p-2};\Z/m)\to E_2^{3,2p-2}\to 0.
    \end{equation}
    Combining \eqref{eq:2p+1_mod_m} and \eqref{eq:SES_2p+3}, we have 
    \begin{equation}\label{eq:Einfty_2p+3}
        E_{\infty}^{2p+1,0}\cong\Z/p\textrm{ generated by }u_{m,p-1}.
    \end{equation}
    Combining \eqref{eq:modEinfty2p+1} and \eqref{eq:Einfty_2p+3}, we have 
    \begin{equation*}
        \rho_m(\kappa_{2p+1}) \equiv  q_{p-1}[\bxxi{m}{r}\om{r}{2}{r} + \sum_{i=1}^{r-1}\bom{r}{1}{i}\om{r}{i+1}{r} + \bom{r}{1}{r-1}\eeta{m}{r}]\pmod{u_{m,p-1}}
    \end{equation*}
    and we conclude.
\end{proof}

Recall the injective homomorphism  $\phi\colon\Z/n\to\Z/m$ defined by $\phi(1) = mn^{-1}$. We abuse notations to let $\phi$ also denote the induced map $\K{n}\to\K{m}$.
\begin{prop}\label{pro:lifting}
    Suppose that we have $1 \leq r \leq p-1$, $n^r\mid m$. Then there is a map $g_{2r}\colon \K{n}\to X_{2r}$ such that $f_{2r}g_{2r}\simeq\phi$. Furthermore, there are integers $l_r$ satisfying  $l_r\equiv q_r\pmod{ n}$, where $q_r$ is defined in Proposition~\ref{pro:Massey_in_X} such that  
    \begin{equation}\label{eq:lifting}
        g_{2r}^*(\kappa_{2r+3}) = \frac{l_rm}{n^r}\zzeta{n}{r+1}\in H^{2r+3}(\K{n},\Z).
    \end{equation}
\end{prop}
\begin{proof}
     We proceed by induction on $r$. For $r=1$, we have $X_{2r}=X_2=\K{m}$, and $f_{2r} =f_2$ is the identity. Therefore, we take $g_2 = \phi$. By Lemma~\ref{lem:cyclic}, (\ref{ite:cyclic1}), the k-invariant $\kappa_5$ generates the group $H^5(\K{m};\Z)$. This means, by Proposition~\ref{pro:Kn_integral}, that there is an integer $q$, coprime to $p$, satisfying $\kappa_5 = q\zzeta{m}{2}$. By Proposition~\ref{pro:Kn_to_Km}, we have 
    \[\phi^*(\zzeta{m}{2}) = \frac{m}{n}\zzeta{n}{2}\]
    and the case $r = 1$ is settled.

    Suppose $2\leq r\leq p-2$ and $n^r\mid m$. By the induction hypothesis, we have 
    \[g_{2r-2}\colon \K{n}\to X_{2r-2}\]
    and $l_{r-1}$ satisfying $l_{r-1}\equiv q_{r-1}\pmod{n}$, and
    \[g_{2r-2}^*(\kappa_{2r+1}) = \frac{l_rm}{n^{r-1}}\zzeta{n}{r} \in H^{2r+1}(\K{n};\Z).\]
    Since $n^r\mid m$ and $n\zzeta{n}{r-1} = 0$, we have $g_{2r-2}^*(\kappa_{2r-1}) =0$. Therefore, we have a map $g_{2r}\colon \K{n}\to X_{2r}$ satisfying $f_{2r}g_{2r}\cong\phi$.

    The homotopy $f_{2r}g_{2r}\cong\phi$ indicates
    \begin{equation}\label{eq:L_to_M}
        \begin{split}
            & g_{2r}^*(L_{m,r}) = g_{2r}^*(\langle x_{m,r}, \underbrace{y_{m,r},\cdots,y_{m,r}}_{r\textrm{ tuples}}\rangle) \subseteq \phi^*(\langle x_m, \underbrace{y_m,\cdots,y_m}_{r\textrm{ tuples}}\rangle)\\
            \subseteq &\langle x_n, \underbrace{y_n,\cdots,y_n}_{r\textrm{ tuples}}\rangle = M_r\subset H^{2r+3}(\K{n};\Z/m). 
        \end{split}        
    \end{equation}
    By Proposition~\ref{pro:xyMassey} and \eqref{eq:L_to_M}, we have 
    \begin{equation}\label{eq:g_of_L}
        g_{2r}^*(L_{m,r})\subseteq M_r =\{\frac{km}{r!n^r}\yy{n}{r+1}\mid k\equiv 1\pmod{n}\}.
    \end{equation}
    By Proposition~\ref{pro:Massey_in_X}, we have
    \[\rho_m(q_r^{-1}\kappa_{2r+3})\in L_{m,r},\]
    and since $\yy{n}{r+1} = \rho_m(\zzeta{n}{r+1})$, we have
    \begin{equation}\label{eq:rho_kappa_rho_z}
        \rho_m(g_{2r}^*(\kappa_{2r+3})) = q_r\rho_m(\frac{km}{r!n^r}\zzeta{n}{r+1}),\ k\equiv 1\pmod{n}.
    \end{equation}
    Since 
    \[\rho_m\colon H^{2r+3}(\K{n};\Z)\to H^{2r+3}(\K{n};\Z/m)\]
    is an isomorphism, we may remove the $\rho_m$ on both sides of \eqref{eq:rho_kappa_rho_z}, and obtain \eqref{eq:lifting}. The inductive proof for $1\leq r\leq p-2$ is concluded.

    For $r = p-1$, the existence of $g_{2p-2}$ such that $f_{2p-2}g_{2p-2} \simeq \phi$ follows from \eqref{eq:lifting} for $r = p-2$. The computation \eqref{eq:L_to_M} hold for $r = p-1$:
    \begin{equation}\label{eq:L_to_M_p}
        g_{2p-2}(L_{m,p})\subset M_{m,p} = \{\frac{km}{(p-1)!n^{p-1}}\yy{n}{p}\mid k\equiv1\pmod{n} \}.
    \end{equation}
    We recall the following facts.
    \begin{itemize}
        \item By Proposition~\ref{pro:generators_for_HK}, we have 
        \[H^{2p+1}(\K{n};\Z/m)\cong \Z/(pn)\oplus\Z/p\]
        where the first summand is generated by $\yy{n}{p} = \rho_m(\zzeta{n}{p})$ and the second one is generated by $u_n$.
        \item The image of $\rho_m\colon H^{2p+1}(\K{n};\Z)\to H^{2p+1}(\K{n};\Z/m)$ is generated by $\rho_m(\zzeta{n}{p}) = \yy{n}{p}$.
        \item By Lemma~\ref{lem:2p+1_X}, we have
        \begin{equation}\label{eq:f_of_w}
            f_{2p-2}^*(u_{m,p-1}) = u_n.
        \end{equation}
        \item By Proposition~\ref{pro:Massey_in_X}, we have $l\in\Z$, $\mu\in L_{m,p}$ such that
        \begin{equation}\label{eq:kappa2p+1}
            \rho_m(\kappa_{2p+1}) = q_{p-1}\mu + lu_{m,p-1}
        \end{equation}
    \end{itemize}
    Combining \eqref{eq:L_to_M_p} and \eqref{eq:kappa2p+1}, we have
    \begin{equation}\label{eq:rho_g}
        \rho_mg_{2p-2}^*(\kappa_{2p+1}) = \frac{q_{p-1}k m}{(p-1)!n^{p-1}}\yy{n}{p} + lu_{p-1},\ k\equiv 1\pmod{m}.
    \end{equation}
    Since the $\Zp$-module
\[\opn{Im}\{\rho_m\colon H^{2p+1}(\K{n};\Zp)\to H^{2p+1}(\K{n},\Z/m)\}\]
is generated by $\yy{n}{p}$, we have $l=0$, and the proof is completed.

\end{proof}
Theorem~\ref{thm:cw}, particularly part (2), follows from Proposition~\ref{pro:lifting}, as we shall see.
\begin{proof}[Proof of Theorem~\ref{thm:cw}]
    Take $X = \K{n}$, and $\alpha\in\Br(X)$ the restriction of $\zeta_{n}$. Then we have
    \begin{equation*}
        E_2^{2i+1,-2i}\cong H^{2i+1}(\K{n};\Zp)\cong
        \begin{cases}
            \Z/(n),\ i < p,\\
            \Z/(pn),\ i = p.
        \end{cases}
    \end{equation*}
   Notice that a generator of $E_{d+1}^{0,0}$ divides the product of the cardinalates of all $E_2^{2i+1,-2i}$, for $2i+1\leq d$. Part (1) of the theorem then follows from Proposition~\ref{pro:per_ind_ahss}. 

   For part (2), take $X = \sk{d}{n}$, the $d$skeleton of $\K{n}$, and $\alpha\in\Br(X)$ the restriction of $\zeta_{n}$. Then Proposition~\ref{pro:lifting} implies that in the diagram \ref{eq:lift_diag_general},
    the dashed arrow exists only if 
    \begin{itemize}
        \item $r < p$, $d\leq 2r+2$ and $n^r\mid m$, or
        \item $d = 2p+2$ and $pn^p\mid m$.
    \end{itemize}
    In particular, this shows that
    \begin{itemize}
        \item for $d < 2p+1$,  a principal $PU(m)$-bundle over $X$ of which the associated Brauer class is $\alpha$exists only if $n^{[(d+1)/2]-1}\mid m$, and
        \item for $d = 2p+1$ or $d = 2p+2$, a principal $PU(m)$-bundle over $X$  of which the associated Brauer class is $\alpha$ exists only if $pn^{[(d+1)/2]-1}\mid m$.
    \end{itemize}

    This proves 
    \begin{equation*}
        \begin{cases}
            \per(\alpha)^{[(d+1)/2]-1}\mid\ind(\alpha),\ d < 2p + 1,\\
            p\per(\alpha)^{[(d+1)/2]-1}\mid \ind(\alpha) ,\ d = 2p +1 \textrm{ or }2p+2.
        \end{cases}
    \end{equation*}
    Combining (1) of Theorem~\ref{thm:cw}, this proves (2) of Theorem~\ref{thm:cw}. 
\end{proof}



\section{Hirzebruch's formula for Steenrod operations}\label{sec:Steenrod}
Let $M$ be a closed oriented differentiable manifold of dimension $2p+2$. Consider the Steenrod power operation
\[\pp^1\colon H^4(M;\F_p)\to H^{2p+2}(M;\F_2).\]
By \cite{hirzebruch1953steenrod}, there is a class $v\in H^{2p-2}(M;\Z)$, which is a polynomial in the Pontrjagin classes of $M$, satisfying 
\begin{equation}\label{eq:Hirzebruch}
    \pp^1(t) = \rho_p(v)t,\ t\in H^4(M;\F_p).
\end{equation}

\begin{remark}\label{eq:v}
    Let $\mathfrak{p}_i$ denote the $i$th Pontrjagin class of $M$. A routine computation shows
    \[\rho_p(v) = \lambda\mathfrak{p}_{(p-1)/2} + \textrm{(other monomials in $\mathfrak{p}_i$ of degree $(p-1)/2$),}\]
    where 
    \begin{equation*}
        \lambda = 
        \begin{cases}
            1,\ p = 3,\\
            3^{(1-p)/2} = \pm 1,\ p > 3.
        \end{cases}
    \end{equation*}
\end{remark}

\section{Linking Pairing on closed oriented manifolds}\label{sec:linking}
We have a bilinear form $\Q/\Z\times\Z\to\Q/\Z$, $(a,b\mapsto ab)$, which gives rise to a cup product
\begin{equation}\label{eq:cup}
    H^s(X;\Q/\Z)\times H^t(X;\Z)\to H^{s+t}(X;\Q/\Z)
\end{equation}
for any topological space $X$. 

Let $TH^k(X)$ denote the torsion subgroup of $H^K(X;\Z)$.
Let $M$ be a closed oriented manifold of dimension $d$, and let $[M]\in H_d(M;\Q/\Z)$ be the fundamental class. Then we have a bilinear form 
\begin{equation}\label{eq:bilinear}
    \td{\lambda}\colon H^k(M;\Q/\Z)\times H^{d-k}(M;\Z)\to \Q/\Z,\ (x,y)\mapsto \langle xy, [M]\rangle
\end{equation}
where the product $xy\in H^d(M;\Q/\Z)$ is defined in the sense of \eqref{eq:cup}.

Let $\delta\colon H^k(M;\Q/\Z)\to H^{k+1}(M;\Z)$ denote the Bockstein homomorphism associated to the short exact sequence
\[0\to\Z\to\Q\to\Q/\Z\to 0.\]
The standard long exact sequence argument shows that the image of $\delta$ is $TH^{k+1}(M)$. Therefore, for $x\in TH^{k+1}(M)$, there is a $\td{x}\in H^k(M;\Q/\Z)$ such that $\delta(\td{x}) = x$. 

The \textit{linking paring} for $M$ is defined by 
\begin{equation}\label{eq:linking}
   \lambda\colon TH^{k+1}(M)\times TH^{d-k}(M)\to \Q/\Z,\ (x,y)\mapsto \td{\lambda}(\td{x},y) = \langle \td{x}y,[M]\rangle.
\end{equation}
By Poincar{\'e} duality, the value of $\lambda(x,y)$ is independent of the choice of $\td{x}$. 

For any group $G$, let $G^{\wedge} = \opn{Hom}(G,\Q/\Z)$. Then a homomorphism
\[B\colon G\times H\to \Q/\Z\]
induces homomorphisms $B_l\colon G\to H^{\wedge}$ and $B_r\colon H\to G^{\wedge}$. We say that $B$ is a \textit{perfect pairing} if both $B_l$ and $B_r$ are isomorphisms.  

\begin{lem}[\cite{crowley2020topological}, Lamma~3.4]\label{lem:perfect_pairing}
    The linking pairing is a perfect pairing.
\end{lem}

For the reminder of this section, we fix a prime number $p$ and let $d = 2p+2$. We have a homomorphism of abelian groups
\[\phi\colon \F_p\to\Q/\Z, \ k\to k/p.\]

Suppose $x\in H^2(M;\F_p)$, we define a homomorphism
\[\lambda_x\colon TH^2(M)\to\Q/\Z,\ y\mapsto \phi(\langle x\rho_p(y)^p, [M]\rangle).\]

\section{On multi-linear forms}\label{sec:bilinear}
Let $V$ be a finite dimensional vector space over a field $\F$. 
\begin{lem}\label{lem:Vdecomp}
   Let $\theta\colon V^{\times n}\to\F$ be a multi-linear form. Then we have direct sum decomposition $V = V_1\oplus V_2$ such that
   \begin{itemize}
       \item For any $\varphi\in V^*$, there are $w_1\cdots,w_{n-1}\in V$ such that for any $v\in V_1$, we have
       \[\varphi(v) = \theta(w_1,\cdots,w_{n-1},v);\]
       \item for any $w_1,\cdots,w_{n-1}\in V$ and $v\in V_2$, we have 
       \[\theta(w_1,\cdots,w_{n-1},v) = 0.\]
   \end{itemize}
\end{lem}
\begin{proof}
    Consider the multi-linear mapping
    \[\mu\colon V^{\times (n-1)}\to V^*,\ (w_1,\cdots,w_{n-1})\mapsto \theta((w_1,\cdots,w_{n-1},-).\]
    Let $V_1'$ be the image of $\mu$, and let $V_2'$ be a linear complement of $V_1'$. Take a dual decomposion $V = V_1\oplus V_2$, which in particular, satisfies $\langle V_1',V_2\rangle = 0$ and $\langle V_2',V_1\rangle=0$. Then the subspaces $V_1, V_2$ are as desired.  
\end{proof}
For the rest of this section, let $\F = \F_p$. 
Let $\theta\colon V^{\times p}\to \F_p$ be a symmetric multi-linear map. Then, we have a linear map $\gamma(\mu)$ defined by
\begin{equation*}
    \gamma(\theta)\colon V\to\F_p,\ v\mapsto \theta(\underbrace{v,\cdots,v}_{\textrm{$p$ tuples}}).
\end{equation*}
\begin{lem}\label{lem:pre_multilinear}
    There are $w_1,\cdots w_{n-1}\in V$, such that for any $v\in V$, we have 
    \[\gamma(\theta)(v) = \theta(w_1,\cdots,w_{p-1},v).\]
\end{lem}
\begin{proof}
    Take $n = p$, and apply Lemma~\ref{lem:Vdecomp} to obtain the decomposition $V = V_1\oplus V_2$. Then we have $w_1,\cdots,w_{p-1}\in V$ such that for $v\in V$ and the decomposition $v = v_1 + v_2$, $v_i\in V_i$, we have
    \begin{equation*}
        \begin{split}
            & \gamma(\theta)(v) = \gamma(\theta)(v_1) + \gamma(\theta)(v_2)\\
            = &\theta(w_1,\cdots,w_{p-1},v_1) + \theta(v_2,v_2,\cdots,v_2)\\
            = &\theta((w_1,\cdots,w_{p-1},v_1)\\
            = &\theta((w_1,\cdots,w_{p-1},v_1) + \theta((w_1,\cdots,w_{p-1},v_2)\\
            =&\theta((w_1,\cdots,w_{p-1},v_1 + v_2) = \theta((w_1,\cdots,w_{p-1},v).
        \end{split}
    \end{equation*}
\end{proof}
Let $V = TH^2(M)/pTH^2(M)\subset H^2(M;\F_p)$, $x\in H^2(M;\F_p)$ and let $\theta$ be defined by 
\begin{equation*}
    \theta(y) = \langle xy^p,[M]\rangle.
\end{equation*}
By Lemma~\ref{lem:pre_multilinear}, we have
\begin{lem}\label{lem:multilinear}
   Given $x\in H^2(M;\F_p)$, there are $x_1,\cdots,x_{p-1}\in TH^2(M)$, such that for any $y\in TH^2(M)$, we have
   \[xy^p = x\rho_p(x_1\cdots x_{p-1})\rho_p(y).\]
\end{lem}

\section{On closed oriented differentiable manifolds of dimension $2p+2$}\label{sec:2p+2}
By \eqref{eq:Hirzebruch}, for $x,y\in H^2(M;\F_p)$, we have 
\begin{equation}\label{eq:vxz}
    \rho_p(v)xz = \pp^1(xz) = x^pz + xz^p.
\end{equation}
If, in addition, we have $z = \rho_p(y)$ for $y\in TH^2(M)$, by Lemma~\ref{lem:multilinear}, we have
$x_1,\cdots,x_{p-1}\in TH^2(M)$ such that
\begin{equation*}
    \rho_p(v)x\rho_p(y) = x^p\rho_p(y) + x\rho_p(x_1\cdots x_{p-1})\rho_p(y),
\end{equation*}
or equivalently,
\begin{equation}\label{eq:vxy}
    [x^p + x\rho_p(x_1\cdots x_{p-1}) - x\rho_p(v)]\rho_p(y) = 0.
\end{equation}
We fix the following notations. Let $\delta_p$ denotes the Bockstein homomorphism associated to the short exact sequence $0\to\Z\xrightarrow{\times n}\Z\to\Z/p\to 0$. Let $\delta_{\Q}$ denote the Bockstein homomorphism associated to the short exact sequence $0\to\Z\to\Q\to \Q/\Z\to 0$. Let $\phi\colon\F_p\to\Q/\Z$ be the homomorphism of abelian groups defined by $k\to k/p$.

\begin{lem}\label{lem:delta}
    \[\delta_p = \delta_{\Q}\phi_*\colon H^*(-;\F_p)\to H^{*+1}(-;\Z).\]    
\end{lem}

\begin{lem}\label{lem:delta_of_product}
    For a topological space $X$, $x\in H^*(X;\F_p)$ and $y\in H^*(X;\Z)$, we have
    \[\delta_p(x_1\rho_p(y)) = \delta_p(x)y.\]  
\end{lem}

\begin{lem}\label{lem:phi}
    For $u\in H^*(M;\F_p)$ and $v\in H^*(X;\Z)$, we have 
    \[\phi(\langle u\rho_p(v),[M]\rangle) = \langle \phi_*(u)v,[M]\rangle.\]
\end{lem}

Applying Lemma~\ref{lem:delta} and Lemma~\ref{lem:delta_of_product} to \eqref{eq:vxy}, we have
\begin{equation}\label{eq:xylinking}
    \begin{split}
        & \delta_{\Q}\phi_*(x^p + x\rho_p(x_1\cdots x_{p-1}) - x\rho_p(v))\\
        = & \delta_p(x^p + x\rho_p(x_1\cdots x_{p-1}) - x\rho_p(v))\\
        = & \delta_p(x^p) + \delta_p(x)x_1\cdots x_{p-1} - \delta_p(x)v.
    \end{split}
\end{equation}
Let $w\in H^2(M;\Z/n)$, where $p\mid n$, and take $x = \rho_{n,p}(w)$. By \eqref{eq:xylinking} and Proposition~\ref{pro:O_n}, we have
\begin{equation*}
    \begin{split}
        & \delta_{\Q}\phi_*(\rho_{n,p}(w)^p + \rho_{n,p}(w)\rho_p(x_1\cdots x_{p-1}) - \rho_{n,p}(w)\rho_p(v)) \\
        =& n\O_p(w) + \delta_p\rho_{n,p}(w)(x_1\cdots x_{p-1} - v).
    \end{split}   
\end{equation*}
A routine computation shows $\delta_p\rho_{n,p} = \frac{n}{p}\delta_n$, and the above becomes
\begin{equation}\label{eq:xylinking1}
    \begin{split}
        & \delta_{\Q}\phi_*(\rho_{n,p}(w)^p + \rho_{n,p}(w)\rho_p(x_1\cdots x_{p-1}) - \rho_{n,p}(w)\rho_p(v)) \\
        =& n\O_p(w) + \frac{n}{p}\delta_n(w)(x_1\cdots x_{p-1} - v).
    \end{split}   
\end{equation}
Applying Lemma~\ref{lem:phi} to \eqref{eq:vxy}, we have
\begin{equation}\label{eq:xylinking2}
    \phi_*(x^p + x\rho_p(x_1\cdots x_{p-1}) - x\rho_p(v))y = 0,\ \forall y\in TH^2(M).
\end{equation}

Recall the linking pairing $\lambda$. Combining \eqref{eq:vxy} and \eqref{eq:xylinking2}, we have
\begin{equation}\label{eq:xylinking3}
    \lambda(n\O_p(w) + \frac{n}{p}\delta_n(w)(x_1\cdots x_{p-1} - v), y) = 0,\ \forall y\in TH^2(M).     
\end{equation}
By Lemma~\ref{lem:perfect_pairing}, this implies
\begin{prop}\label{pro:nO_p}
    Suppose $p\mid n$. There are $x_1,\cdots,x_{p-1}\in TH^2(M)$,  such that for all $w\in H^2(M;\Z/(n))$ we have
    \begin{equation*}
    n\O_p(w) = \frac{n}{p}\delta_n(w)(v - x_1\cdots x_{p-1}).
\end{equation*}
\end{prop}

\begin{proof}[Proof of Theorem~\ref{thm:manifold}]
    It suffices to consider the case where $M$ is of dimension $2p+2$. Let $\alpha\in\Br(M)$ with $\per(\alpha) = n$, and let $w\in H^2(M;\Z/(n))$ such that $\delta_n(w) = \alpha$.

    Recall that we have 
    \[E_2^{2p+1.-2p}\cong H^{2p+1}(\K{n};\Zp)\cong\Z/np.\]
    By Proposition~\ref{pro:per_ind_ahss}, it suffices to show that the torsion order of $E_{2p+1}^{2p+1.-2p}$ divides $n$. Let $v$ and $x_i$ be as in Proposition~\ref{pro:nO_p}.

    By the product formula, i.e.,  Proposition~\ref{pro:product}, we have 
    \[d_3(x_1\cdots x_p - v) = \widetilde{d}_3((x_1\cdots x_p - v)1 + (x_1\cdots x_p - v)\alpha.\]
    Since we are working $p$-locally, we have $\widetilde{d}_3(y) = 0$ for any even dimensional element in $\widetilde{E}_3$, which follows from the universal case of $K(\Z,2q)$. Therefore, we have
    \[d_3(x_1\cdots x_p - v) =  (x_1\cdots x_p - v)\alpha,\]
    and by Proposition~\ref{pro:nO_p}, we have 
    \[ n\O_p(w) = \frac{n}{p} = \frac{n}{p}d_3(x_1\cdots x_p - v),\]
    and this shows that the torsion order of $E_{2p+1}^{2p+1.-2p}$ divides $n$.
\end{proof}

\section{On the cohomology of $BG(n,m)$}\label{sec:G}
Throughout this section, let $n$ be a power of $p$ and let $n\mid m$. Recall the notation $G(n,m) = SU(m)/\mu_n$. The short exact sequence 
\[1\to\mu_n\to SU(m)\to G(n,m)\to 1\]
gives rise to a homotopy fiber sequence
\[BSU(m)\to BG(n,m)\xrightarrow{\alpha_{n,m}} K(\Z/n,2).\]

In this section we prove
\begin{thm*}[Theorem~\ref{thm:Gnm}] 
    Suppose one of the following is true:
    \begin{enumerate}
        \item $r < p$, $n$ a power of $p$, $n^r\mid m$, or
        \item $r= p$, $n$ a power of $p$, $pn^r\mid m$. 
    \end{enumerate}
    Then the map
    \[\alpha_{n,m}^*\colon H^{\leq 2r+2}(K(\Z/n,2);\Z)\to  H^{\leq 2r+2}(BG(n,m);\Z)\]
    is injective.
\end{thm*}

Then the quotient map $BSU(m)\to BG(n,m)$ is a $\mu_n$-torsor, which is classified by a map
\[\alpha_{n,m}\colon BG(n,m)\to K(\Z/n,2).\]
Let $\iota_{n,m}\colon K(\Z/n,2)\to K(\Z/m,2)$ denote the map induced by the inclusion $\Z/n\hookrightarrow\Z/m$. Let $\beta_m\in H^2(BPU(m);\Z/m)$ be the canonical generator.

\begin{lem}\label{lem:BG_bundle}
   Let $X$ be a CW-complex, and let $\alpha\colon K(\Z/n,2)$ be a map. If there is a map $f\colon X\to BPU(m)$ satisfying $\iota_{n,m}\alpha_{n,m}\simeq\beta_mf$, then $f$ factors, up to homotopy, as
   \[X\xrightarrow{f'}BG(n,m)\to BPU(m),\]
   where the second map is induced by the quotient map.

   If, in addition, $H^2(X;\Z/n)\to H^2(X;\Z/m)$ is injective, then we have $\alpha\simeq\alpha_{n,m}f'$. 
\end{lem}
\begin{proof}
    We have the following homotopy commutative diagram:
    \begin{equation*}
        \begin{tikzcd}
             & & BG(n,m)\arrow[d]\arrow[ddl,"\alpha_{n,m}"description]&\\
             & & BPU(m)\arrow[d,"\beta_m"description]\arrow[dr,"\beta'"description] &\\
            X\arrow[uurr,dashed,"f'"description]\arrow[urr,"f"description]\arrow[r,"\alpha"description]&K(\Z/n,2)\arrow[r,"\iota_{n,m}"]&K(\Z/m,2)\arrow[r]&K(\Z/\frac{m}{n},2)
        \end{tikzcd}
    \end{equation*}
    The existence of the dashed arrow $f'$ follows from the fact that 
    \[K(\Z/n,2)\xrightarrow{\iota_{n,m}}K(\Z/m,2)\to K(\Z/\frac{m}{n},2)\]
    and 
    \[BG(n,m)\to BPU(m)\xrightarrow{\beta'}K(\Z/\frac{m}{n},2)\]
    are homotopy fiber sequences. The injectivity statement is a routine computation.
\end{proof}
\begin{proof}[Proof of Theorem~\ref{thm:Gnm}]
    For $r\leq p-1$, $n^r\mid m$, or $r = p$, $pn^p\mid m$, let $X$ be the $(2r+2)$th skeleton of $K(\Z/n,2)$. Then Lemma~\ref{lem:BG_bundle} yields a decomposition of the inclusion $X\hookrightarrow K(\Z/n,2)$ as 
\[X\xrightarrow{f'} BG(n,m)\xrightarrow{\alpha_{n,m}}K(\Z/n,2).\]
The injectivity follows upon applying cohomology to this decomposition.
\end{proof}

Let $\alpha_{n,m}'\in H^2(BG(n,m);\Z/m)$ be the image of $\alpha_{n,m}$ under the canonical homomorphism $H^2(BG(n,m);\Z/n)\to H^2(BG(n,m);\Z/m)$. Let ${\alpha}_{n,m}''\in H^3(BG(n,m);\Z/m)$ be the image of $\alpha_{n,m}$ under the Bockstein homomorphism associated to
\[0\to\Z/m\to\Z/mn\to\Z/n\to 0.\]
\begin{cor}\label{cor:Massey_in_BG}
   For $r\leq p-1$, $n^r\mid m$, or $r = p$, $pn^p\mid m$, the Massey product
   \[N_{n,m,r} = \langle\alpha_{n,m}',\underbrace{\alpha''_{n,m},\cdots,\alpha''_{n,m}}_{r\textrm{ tuples}}\rangle\]
   is well defined, and does not contain $0$. 
\end{cor}

\bibliographystyle{abbrv}
\bibliography{references.bib}
\end{document}